\documentclass[11pt,a4paper, parskip=half]{scrartcl} 
\usepackage[utf8]{inputenc}
\usepackage[english]{babel}
\usepackage[T1]{fontenc}
\usepackage{lmodern}
\usepackage[left=3cm,right=3cm,top=3cm,bottom=3cm]{geometry}
\usepackage[runin]{abstract}
    \abslabeldelim{.}

\usepackage{todonotes}

\usepackage{amsmath, amsfonts, amssymb, amsthm, amstext}
\usepackage{mathtools}
\usepackage{mathrsfs}
\usepackage{latexsym}
\usepackage{mathdots}
\usepackage{cases}

\usepackage{graphicx}
\usepackage{subcaption}
\usepackage{algorithm}

 at11pt
 at9pt

\usepackage{enumitem}

\usepackage{color}
\usepackage{tikz}
\usepackage{bm}
\usepackage[numbers,sort]{natbib}

\usepackage{hyperref} 
\usepackage[capitalise]{cleveref}
	
\newcounter{thm}
\numberwithin{thm}{section}
\numberwithin{equation}{section}

	\newtheoremstyle{myplain}		
			{}			
			{}			
			{\itshape}				
			{}				
			{\sffamily\bfseries}				
			{.}		
			{ }				
			{\thmname{#1}\thmnumber{ #2}\textnormal{\textsf{\thmnote{ (#3)}}}}			
    \newtheoremstyle{mybreak}
            {}{}{}{}{\sffamily\bfseries}{.}{\newline}
            {\thmname{#1}\thmnumber{ #2}\textnormal{\textsf{\thmnote{ (#3)}}}}
	\newtheoremstyle{mydef}
			{}{}{}{}{\sffamily\bfseries}{.}{ }
			{\thmname{#1}\thmnumber{ #2}}
	\newtheoremstyle{myrem}
			{}{}{}{}{\sffamily\itshape}{.}{ }
			{\thmname{#1}\thmnumber{ #2}}

\theoremstyle{myplain}

\theoremstyle{mybreak}
\theoremstyle{mydef}
	
	\newtheorem{remark}[thm]{Remark}
\theoremstyle{mydef}
	\newtheorem{example}[thm]{Example}

	\newcommand{\cc}{\mathbb{C}}
		\newcommand{\zz}{\mathbb{Z}}
		\newcommand{\nn}{\mathbb{N}}
	\newcommand{\rr}{\mathbb{R}}

\DeclareMathOperator{\rank}{rank}
\newcommand{\argmax}{\mathop{\mathrm{argmax}}}

\DeclarePairedDelimiter{\abs}{\lvert}{\rvert}

\allowdisplaybreaks

\def\sumprime_#1^#2{
    \setbox0=\hbox{$\scriptstyle{#1}$}
    \setbox1=\hbox{$\scriptstyle{#2}$}
    \setbox2=\hbox{$\displaystyle{\sum}$}
    \setbox4=\hbox{${}^\prime\mathsurround=0pt$}
    \dimen0=.5\wd0 \advance\dimen0 by-.5\wd2
    \ifdim\dimen0>0pt
        \ifdim\dimen0>\wd4 \kern\wd4
        \else\kern\dimen0
        \ifdim\dimen1>\wd4 \kern\wd4
        \else\kern\dimen1
    \fi\fi\fi
\mathop{{\sum}^\prime}_{\kern-\wd4 #1}^{\kern-\wd4 #2}
}
\title{\Large Parameter estimation for multivariate exponential sums via iterative rational approximation}
\author{Nadiia Derevianko\footnote{TUM School of CIT, 
Department of Computer Science,
Boltzmannstrasse 3,
85748 Garching b. München,
Germany,  nadiia.derevianko@tum.de} \footnote{Corresponding author},   \  Lennart Aljoscha Hübner \footnote{KU Leuven, Department of Mathematics, Celestijnenlaan 200b - box 2400, 
3001 Leuven, Belgium, lennart.huebner@kuleuven.be}}
\date{\today}

\begin{document}
	\let\oldproofname=\proofname
	\renewcommand{\proofname}{\itshape\sffamily{\oldproofname}}

\maketitle

\begin{abstract}
We present two new methods for multivariate exponential analysis. In  \cite{DP21},  we developed a new  algorithm for reconstruction of univariate  exponential sums by exploiting the rational structure of their Fourier coefficients and reconstructing   this rational structure with the  AAA (``adaptive Antoulas–Anderson'')  method for rational approximation \cite{AAA}. In this paper, we extend these ideas  to the multivariate setting. Similarly as in univariate case, the Fourier coefficients of multivariate exponential sums have a rational structure and the multivariate exponential recovery problem can be reformulated as multivariate rational interpolation problem. We develop two approaches to solve this special multivariate rational interpolation problem by reducing it to the several univariate ones, which are then solved again via the univariate AAA method. Our first approach is based on using indices of the Fourier coefficients chosen from some sparse grid, which ensures efficient reconstruction using a respectively small amount of input data. The second approach is based on using the full grid of indices of the Fourier coefficients and relies on the idea of recursive dimension reduction. We demonstrate performance of our methods with  several numerical examples. 
\\[1ex]
\end{abstract}

\section{Introduction }

In the paper,  we develop two  new methods for multivariate exponential analysis,  using  connections to rational functions.
Let the dimension $d\in \nn$  be given.  We define a \emph{$d$-variate  exponential sum} of order $M \in \nn$ by 
\begin{equation}\label{mul}
f(\boldsymbol{t})\coloneqq\sum\limits_{j=1}^{M} \gamma_j \mathrm{e}^{\langle \boldsymbol{\lambda}_j, \boldsymbol{t} \rangle},   \  \    \boldsymbol{t}\coloneqq(t_1,\ldots,t_d) \in \rr^{d},
\end{equation}
with complex coefficients $\gamma_j \neq 0$,  distinct frequency vectors  $\boldsymbol{\lambda}_j=(\lambda_{j1},\ldots, \lambda_{j d}) \in \cc^{d}$ and 
 $\langle \boldsymbol{\lambda}_j , \boldsymbol{t} \rangle\coloneqq \sum_{\ell=1}^{d} \lambda_{j\ell} t_\ell$.  Our goal is to reconstruct $M$, 
 the vectors $\boldsymbol{\lambda}_j$,  $j=1,\ldots,M$,  i.e each individual  $\lambda_{j\ell}$,  $j=1,\ldots,M$,  $\ell=1,\ldots,d$,  and the coefficients $\gamma_j$,  $j=1,\ldots,M$,  employing as input data a finite set of  their Fourier coefficients defined for some $P>0$ by
\begin{equation}\label{fcm}
c_{\boldsymbol{k}}(f)\coloneqq\frac{1}{P^d}\int_{[0,P]^{d}} f (\boldsymbol{t}) \, \mathrm{e}^{-\frac{2 \pi \mathrm{i}}{P} \langle \boldsymbol{k},\boldsymbol{t}\rangle} \,  \mathrm{d} \boldsymbol{t}, \quad \boldsymbol{k} \in  \zz^d.
\end{equation}

Multivariate exponential sums have a wide range of applications.  They are used in signal processing \cite{ SD2007,  YFS2006},  fitting
nuclear magnetic resonance spectroscopic data \cite{PLH04},    fast algorithms \cite{G2022},  etc. The connection between multivariate exponential sums and PDEs were presented in the recent paper \cite{sturm2023}. The authors study a  very general case of PDEs (namely, homogeneous linear partial differential equations with constant coefficients) that describe vector-valued functions from $d$-space to $n$-space (in our case $n=1$) under the requirements that the space of  solutions contains the exponential functions  
$
q(t_1,\ldots,t_d) \mathrm{e}^{\lambda_1 t_1+\ldots +\lambda_d t_d},
$
where $q$ is a multivariate algebraic polynomial.

As alternative to the Prony-like methods for univariate exponential recovery (see for example \cite{RT89, PT13, L20, B23}),  
in \cite{DP21} we developed a new method for reconstruction of exponential sums, which is based on exploiting the rational structure of their Fourier coefficients and reconstruction of  this rational structure using  the well-known AAA method for rational approximation \cite{AAA}. Applying the AAA routine in the frequency domain ensures numerical stability and small complexity of the method.  Our goal in this paper is to extend results from \cite{DP21}
 to the multivariate setting.

  Our research is closely related to the  alternative strategy for rational recovery, which relies on the exponential structure of the Fourier transform of rational functions and its reconstruction using Prony-like methods in the frequency domain.  This technique was employed in  \cite{DB13,  WDT21, Y22, D24} for the recovery of rational functions.   Note that  in \cite{D24}, the reconstruction from noisy data is not limited to numerical observations. To explain the reconstruction results achieved using the proposed in \cite{D24} method, the sensitivity analysis of poles was presented.

Multivariate Prony-like methods where considered  in several publications (see, for example, \cite{DI17, KP2016, PV2020, S2018, M25}).   
Of particular interest to us are the following two methods: a sparse approximate Prony method (SAPM)  \cite{PT13mult} and a method for multivariate exponential recovery from the minimal number of samples  \cite{CL2018}.  Both of them  propose recovery using relatively small or minimal amount of given data and  are simple to implement.  In \cite{PT13mult},  the authors developed SAPM based on using $\mathcal{O}(d N)$, $N\geq M$ samples from $\zz^d$ located  on  few straight lines.  
The method proposed in \cite{CL2018}  employs  the minimal number of samples
 $(d+1)M$ and  a 1-dimensional Prony technique.
 

  Straightforward generalization of the results from \cite{DP21} will bring us to the multivariate rational interpolation problem (\ref{intd}), which does not allow
the direct computation of
the poles (see Remark \ref{rmrf} for further explanations). Therefore, inspired by the ideas of \cite{PT13mult} and \cite{CL2018}, we reduce the multivariate rational interpolation problem to several univariate ones and solve them using the AAA method  or the Loewner pencil method for univariate  rational approximation. We develop two different approaches to reduce the multivariate rational interpolation problem (\ref{intd}) to several univariate ones. Let $N \in \nn$ be given  and $N\geq M$. We consider the Fourier coefficients $c_{\boldsymbol{k}}(f)$ in (\ref{fcm}) with indices $\boldsymbol{k} \in [-N,N]^d \cap \zz^d$. Our first approach is based on using only  $\mathcal{O}(d N)$ Fourier coefficients $c_{\boldsymbol{k}}(f)$. Indices $\boldsymbol{k}$ of the corresponding Fourier coefficients are chosen from $2d-1$ straight lines and create a \textit{sparse grid} in the $d$-dimentional cube $[-N,N]^d \cap \zz^d$. We use $d$ lines to  separately reconstruct sets of frequencies with respect to each dimension by applying the AAA method or Loewner pencil approach for univariate rational approximation, i.e. to reconstruct sets $\{\lambda_{j1}\}_{j=1}^M$,..., $\{\lambda_{j d}\}_{j=1}^M$. After that we use $(d-1)$ other lines to create the correct pairs $(\lambda_{j1},\ldots, \lambda_{j d})$, $j=1,...,M$. Our second method is based on using the full grid of input data, i.e. we employ  $\mathcal{O}(N^d)$ Fourier coefficients.  We apply the idea of recursive dimension reduction and reconstruction of frequencies with respect to each dimension by a univariate rational method. The sparse grid approach has an overall computational cost of $\mathcal{O}(d N M^3)$, while the  recursive dimension reduction method requires $\mathcal{O}(N M^2((d-1) M^2 + N^{d-1} ))$ flops. The sparse grid method ensures efficient reconstruction using a respectively small amount of input data, but is restricted to the recovery of exponential sums that satisfy $\lambda_{j m}\neq \lambda_{i m}$ for $j \neq i$ and $m=1,...,d$. The second approach is computationally more expensive and requires the full grid of the Fourier coefficients, but can handle the recovery of any exponential sum with pairwise distinct frequency vectors.


  The rest of the paper is organized as follows. In Section 2, we present  ideas of the  recovery of univariate exponential sums via rational interpolation in the frequency domain, as developed in \cite{DP21}.  In Section 3, we show that similarly to the univariate case, we can reformulate the multivariate parameter estimation problem (\ref{mul}) via the multivariate rational interpolation problem. In Section 4, we describe the main ideas of our methods for the bivariate case and in Section 5, we extend these ideas to any dimension $d>2$. Finally, in Section 6, we present numerical experiments. We finish the paper with the Conclusions section, where we present further modifications of our new methods and plans for the future work.

\section{Recovery of univariate exponential sums via iterative rational approximation } 
\label{recuni}

Let a \textit{univariate exponential sum} of order $M \in \mathbb{N}$ be given by
\begin{equation}\label{unexp}
f(t)= \sum\limits_{j=1}^{M} \gamma_j \mathrm{e}^{\lambda_j t},   \quad t\in \mathbb{R},
\end{equation}
with $\gamma_j,  \lambda_j \in \mathbb{C}$,  $\gamma_j \neq 0$ and $\lambda_j$ be pairwise distinct.  We consider the problem of parameters estimation for the exponential sums (\ref{unexp}) (i.e.  the problem of the recovery of all parameters that determine (\ref{unexp})),  using a finite set of its classical  Fourier coefficients $c_k(f)$. More precisely, we consider the Fourier expansion of $f$ in (\ref{unexp})  in $[0,P]$ for some $P>0$,
$$
f(t)= \sum \limits_{k\in  \mathbb{Z}} c_k(f) \mathrm{e}^{2 \pi \mathrm{i} k t/P}, 
$$
with the Fourier coefficients $c_k(f)= \frac{1}{P}\int\limits_0^{P} f(t)  \, \mathrm{e}^{-2 \pi \mathrm{i} k t/P} \,  \mathrm{d} t $. 
For a single exponent $\phi(t)=\gamma_j \mathrm{e}^{\lambda_j t}$,   we have 
\begin{equation}\label{coef}
c_k(\phi)=\frac{1}{P}\int\limits_0^{P} \gamma_j \mathrm{e}^{\lambda_j t} \mathrm{e}^{-2 \pi \mathrm{i} k t/P} \,  \mathrm{d} t = \frac{\gamma_j}{P}   \int\limits_0^{P}  \mathrm{e}^{(\lambda_j-2 \pi \mathrm{i} k/P ) t} \,  \mathrm{d} t =\begin{cases}
\frac{\gamma_j (\mathrm{e}^{\lambda_j P}-1)}{\lambda_j P-2\pi \mathrm{i} k} ,  &  \lambda_j \neq 2\pi \mathrm{i} k/P, \\
\gamma_j ,    & \lambda_j = 2\pi \mathrm{i} k/P .
\end{cases}
\end{equation}
Using the last formula, for $\lambda_j \neq 2\pi \mathrm{i} k/P$  we get 
\begin{equation}\label{uncoef}
c_k(f)=  \sum\limits_{j=1}^{M} \frac{\gamma_j}{P}  \int\limits_0^{P} \mathrm{e}^{(\lambda_j-2 \pi \mathrm{i} k/P ) t} \,  \mathrm{d} t =\sum\limits_{j=1}^{M} \frac{\gamma_j (\mathrm{e}^{\lambda_j P}-1)}{\lambda_j P-2\pi \mathrm{i} k}=\sum\limits_{j=1}^{M} \frac{\frac{\gamma_j}{2\pi \mathrm{i}}( 1- \mathrm{e}^{\lambda_j P})}{k- \lambda_j P/2\pi \mathrm{i}}.
\end{equation}
Define now a rational function $r(z)$ of type $(M-1,M)$  by its pole-residue representation
\begin{equation}\label{unrat}
r(z)=\sum\limits_{j=1}^{M} \frac{a_j}{z-b_j}.
\end{equation}
  For the choice of parameters $a_j$ and $b_j$ as
$$
a_j=\frac{\gamma_j}{2\pi \mathrm{i}}( 1- \mathrm{e}^{\lambda_j P}),  \quad  b_j= \lambda_j P/2\pi \mathrm{i},
$$
from (\ref{uncoef}) we obtain the following interpolation conditions
\begin{equation}\label{unint1}
r(k)=c_k(f),  \quad k \in \mathbb{Z}. 
\end{equation}
That means that the exponential recovery problem (\ref{unexp}) can be reformulated as the rational interpolation problem (\ref{unint1}). Further, we consider two methods to reconstruct poles $b_j$ of a rational function (\ref{unrat}) using as input information samples $r(k)=c_k(f)$, $k=-N,...,N$, $N\geq M$. When the poles $b_j$ are reconstructed, we can compute coefficients $a_j$ as a least squares solution to the system of linear equations 
$$
\sum\limits_{j=1}^{M} \frac{a_j}{k-b_j}=c_k(f), \quad k=-N,...,N.
$$
Then by $\lambda_j= 2 \pi  \mathrm{i} b_j/P$ and $\gamma_j= \frac{2\pi \mathrm{i} a_j}{1- \mathrm{e}^{\lambda_j P}}$ for $j=1,...,M$  we can reconstruct the frequencies $\lambda_j$ and coefficients $\gamma_j$ of the exponential sum (\ref{unexp}), respectively.


\subsection{Recovery of poles via the AAA method for rational approximation }
\label{ESPIRAwithAAA}

In this subsection,  we describe the main ideas for recovering the poles of a rational function (\ref{unrat}) using the AAA method for rational approximation \cite{AAA}.  We have given sets of sample points $\{-N,...,N  \}$ and the corresponding function values $\{c_{-N}(f),...,c_{N}(f)  \}$, $N\geq M$, that are indices of Fourier coefficients and the corresponding Fourier coefficients, respectively. We want to construct a rational approximation $r(z)$ such that 
$$
r(k)=c_k(f),  \quad k=-N,...,N. 
$$
Consider a partition of the index set $\{-N,...,N  \}=\{k_1,...,k_n \} \cup \{ \tilde{k}_1,...,\tilde{k}_{2N+1-n} \}\coloneqq  S_{n} \cup \Gamma_{n}$. Then $|S_n|=n$ and $|\Gamma_{n}|=2N+1-n$, where by $|\mathcal{X}|$ we denote the number of elements in the set $\mathcal{X}$.
Define a rational function $r$ in its barycentric form 
\begin{equation}\label{unbar}
r(z)=  \frac{p(z)}{q(z)} = \sum\limits_{s=1}^{n}  \frac{c_{k_s}(f)w_s}{z-k_s} \bigg/ \sum\limits_{s=1}^{n} \frac{w_s}{z-k_s},
\end{equation}
where $k_s \in S_n$, $c_{k_s}(f)$ are the corresponding Fourier coefficients, and $w_s$ are weights.
This rational function has type $(n-1,n-1)$ and interpolates $c_{k_s}(f)$, i.e. $r(k_s)= c_{k_s}(f)$ for $s=1,...,n$.

The weights $w_s$ in (\ref{unbar}) are computed to minimize the least squares  error over the remaining sampling points $\tilde{k}_{\ell} \in \Gamma_n$, which is a nonlinear problem.  The AAA method uses linearization of this least squares problem: for $\tilde{k}_{\ell} \in \Gamma_n$, we observe
\begin{align*}
    c_{\tilde{k}_{\ell}}(f) - r(\tilde{k}_{\ell}) &  =\frac{1}{q(\tilde{k}_{\ell})} \left( q(\tilde{k}_{\ell}) c_{\tilde{k}_{\ell}} (f) -p(\tilde{k}_{\ell}) \right) \leadsto q(\tilde{k}_{\ell}) c_{\tilde{k}_{\ell}} (f) -p(\tilde{k}_{\ell}) \\
    &      = \sum\limits_{s=1}^{n} \frac{(c_{\tilde{k}_{\ell}}(f)  -c_{k_{s}}(f) ) w_s}{\tilde{k}_{\ell}-k_s}  =  \boldsymbol{e}_\ell^T L_{2N+1-n,n} \boldsymbol{w},
\end{align*}
where $\boldsymbol{e}_\ell \in \mathbb{C}^{2N+1-n}$ denotes the  $\ell$th canonical column vector,  $\boldsymbol{w}=(w_1,...,w_n)^T \in \mathbb{C}^{n}$ is the weight vector and $L_{2N+1-n,n} \in \mathbb{C}^{(2N+1-n)\times n}$ is a Loewner matrix defined as
\begin{equation}\label{loew}
   L_{2N+1-n,n}=  \left( \frac{c_{\tilde{k}_{\ell}}(f)  -c_{k_{s}} (f)}{\tilde{k}_{\ell}-k_s}   \right)_{\ell=1, \ s=1}^{2N+1-n, \ n }.
\end{equation}
Then we compute the weight vector  $\boldsymbol{w}$ as a solution of the linear least squares problem 
\begin{equation}\label{lls}
\min\limits_{\| \boldsymbol{w} \|_2=1} \| L_{2N+1-n,n} \boldsymbol{w} \|_2.
\end{equation}
The AAA algorithm works iteratively with respect to type $(n-1,n-1)$ of the rational approximation (\ref{unbar}), and at each iteration step,  a greedy search is used to select the next point to be added to the set $S_n$. At the step $n$, we have constructed the rational approximant (\ref{unbar}). The next  point $k_{n+1}$ to the set $S_n$ is chosen such that 
$$
k_{n+1}= \argmax \limits_{\ell=1,...,2N+1-n} |c_{\tilde{k}_{\ell}}(f) - r(\tilde{k}_{\ell})|. 
$$
Then this point $k_{n+1}$ is deleted from the set $\Gamma_n$. AAA terminates when either a prespecified error tolerance or the wanted order is achieved.  

To reconstruct a rational function of type $(M-1,M)$ we need $M+1$ iterations  of the AAA algorithm.  After $M+1$ steps we get a rational function 
\begin{equation}\label{finbar}
r(z)= \sum\limits_{s=1}^{M+1}  \frac{c_{k_s}(f)w_s}{z-k_s}  \bigg/ \sum\limits_{s=1}^{M+1} \frac{w_s}{z-k_s},
\end{equation}
that coincides (in case of exact data) with $r$ as in (\ref{unrat}) (see \cite{DP21} and \cite{DPP21} for details). Therefore, the poles $b_j$, $j=1,...,M$ can be reconstructed as poles of (\ref{finbar}). To compute $b_j$, $j=1,...,M$ we just have to solve the following  $(M+2)\times (M+2)$ generalized eigenvalue problem
\begin{equation}\label{eig} \left( \begin{array}{ccccc}
0 & w_{{1}} & w_{{2}} & \ldots & w_{{M+1}} \\
1 & k_{1} &   &     & \\
1 & & k_{2} & & \\
\vdots & & & \ddots & \\
1 & & & & k_{M+1} \end{array} \right) \, {\mathbf v}_{z}= z \left( \begin{array}{ccccc}
0 & & & & \\
& 1 & & & \\
& & 1 & & \\
& & & \ddots & \\
& & & & 1 \end{array} \right) \, {\mathbf v}_{z}.
\end{equation}
Two eigenvalues of this generalized eigenvalue problem are infinite and the other $M$ eigenvalues are the wanted poles  $b_j$ of (\ref{finbar}) (see for example \cite{AAA} for  details). 

\begin{remark}
\label{rem21}
 Note that we can apply the approach described above only for frequencies  $\lambda_j \neq 2\pi \mathrm{i} k/P$, $k \in \mathbb{Z}$. When $\lambda_j = 2\pi \mathrm{i} k/P$ for some $k \in \mathbb{Z}$, then the Fourier coefficients  $c_k(f)$ loose the rational structure (see (\ref{coef})) and can not be reconstructed via a method for rational approximation. Therefore, frequencies of the form $\lambda_j = 2\pi \mathrm{i} k/P$ has to be reconstructed in a post-processing step by comparing given Fourier coefficients $c_k(f)$ and function values $r(k)$ of the computed rational approximation (\ref{finbar}). We refer the interested reader to  \cite{DP21} for further  details.
\end{remark}

\subsection{Recovery of poles via Loewner pencil method}
\label{mpa}

In this subsection, we introduce an approach based on Loewner pencil for computation of poles of a rational function (\ref{unrat}). Note that in \cite{DP21} we have only considered the computation of the poles $b_j$ of (\ref{unrat}) using the AAA method.  

Let again $N \geq M$.   Consider a partition of the index set $\{-N,...,N \}=\{k_1,...,k_M \} \cup \{ \tilde{k}_1,...,\tilde{k}_{2N+1-M} \}\coloneqq  S_{M} \cup \Gamma_{M}$. Then $|S_M|=M$ and $|\Gamma_{M}|=2N+1-M$. Further,  we create two Loewner matrices
\begin{align}
L_{2N+1-M,M}(t)  =\left(  \frac{\tilde{k}_{\ell}^t c_{\tilde{k}_{\ell}}(f)  -k_s^t c_{k_{s}} (f)}{\tilde{k}_{\ell}-k_s}   \right)_{\ell=1, \ s=1}^{2N+1-M, \ M }, \quad t=0,1. \notag  
\end{align}
There are different methods for selecting the points for the sets $S_M$ and $\Gamma_M$. We refer the interested reader to \cite{ALI17, KGA}. In this paper we consider the case when  the set $S_M$ is computed in a greedy fashion by the AAA algorithm.  Then the Loewner matrix $L_{2N+1-M,M}(0)$ is exactly the Loewner matrix created by the AAA algorithm after $M$ iterations and  $L_{2N+1-M,M}(1)$ is the so-called shifted Loewner matrix. 


As we will see shortly, the  Loewner matrices $L_{2N+1-M,M}(t)$, $t=0,1$ offer favourable factorizations.
First, we define the Cauchy matrix $C_{\alpha,\beta}\coloneqq\left( \frac{1}{\alpha_j-\beta_k} \right)_{j=1,k=1}^{n,m} \in \mathbb{C}^{n\times m}$ for two vectors $\alpha=(\alpha_1,\ldots,\alpha_n)$ and $\beta=(\beta_1,\ldots,\beta_m)$.
Then we set $\tilde{\boldsymbol{k}}=(\tilde{k}_{\ell})_{\ell=1}^{2N+1-M}$, $\boldsymbol{k}=(k_s)_{s=1}^{M} $ and $\boldsymbol{b}=(b_j)_{j=1}^{M}$, and obtain
\begin{align}
L_{2N+1-M,M}(0) & = \left( \frac{\sum\limits_{j=1}^{M} \frac{a_j}{\tilde{k}_{\ell}-b_j} -\sum\limits_{j=1}^{M} \frac{a_j}{k_s-b_j}}{\tilde{k}_{\ell}-k_s}   \right)_{\ell=1, s=1}^{2N+1-M,  M }  = \left( \sum\limits_{j=1}^{M} \frac{-a_j}{(\tilde{k}_{\ell}-b_j)(k_s - b_j)}  \right)_{\ell=1, s=1}^{2N+1-M,  M }  \notag \\
&= - \left( \frac{1}{\tilde{k}_{\ell}-b_j}\right)_{\ell=1,  s=1}^{2N+1-M,  M } \,  \mathrm{diag} (a_j)_{j=1}^{M} \,  \left(  \left( \frac{1}{k_s-b_j}\right)_{s=1, \ j=1}^{M, \ M } \right)^{T} \notag \\
& = - C_{\tilde{\boldsymbol{k}}, \boldsymbol{b}} \,  \mathrm{diag} (a_j)_{j=1}^{M} \,   C^{T}_{\boldsymbol{k},  \boldsymbol{b}}.    \label{fac1} 
\end{align}
In a similar manner we get
\begin{align}
L_{2N+1-M,M}(1)  = - C_{\tilde{\boldsymbol{k}}, \boldsymbol{b}} \,  \mathrm{diag} (a_j)_{j=1}^{M} \,  \mathrm{diag} (b_j)_{j=1}^{M} \,   C^{T}_{\boldsymbol{k},  \boldsymbol{b}}.   \label{fac2} 
\end{align}
From (\ref{fac1}) and (\ref{fac2}),  we derive
$$
z L_{2N+1-M,M}(0)  - L_{2N+1-M,M}(1) =  - C_{\tilde{\boldsymbol{k}}, \boldsymbol{b}}  \,  \mathrm{diag} (a_j)_{j=1}^{M} \,   \mathrm{diag} (z-b_j)_{j=1}^{M} \,    C^{T}_{\boldsymbol{k},  \boldsymbol{b}},
$$
which means that the poles $b_j$,  $j=1,\ldots,M$ can be reconstructed as eigenvalues of the Loewner pencil 
\begin{equation}\label{lp}
z L_{2N+1-M,M}(0)  - L_{2N+1-M,M}(1) .
\end{equation}

From (\ref{fac1}) and (\ref{fac2}) follows that both Loewner matrices  $L_{2N+1-M,M}(0)$ as well as $L_{2N+1-M,M}(1)$ have rank $M$. It is easy to see that for any partition $\{-N,...,N  \}=  S \cup \Gamma$ such that $|S|\geq M$ and $|\Gamma|\geq M$, the corresponding Loewner matrices $L_{|\Gamma|,|S|}(0)$ and $L_{|\Gamma|,|S|}(1)$ have ranks $M$. This property does not depend on the partition, $S$ and $\Gamma$, as long they both contain at least $M$ elements. We refer to \cite{ALI17,KGA} regarding this and further properties of Loewner matrices and their connection to rational approximation (see  also \cite{DPP21, DPR23}). 

From Subsection \ref{ESPIRAwithAAA} we know that the AAA routine terminates after $M+1$ iteration steps and creates at the last step a Loewner matrix 
$
   L_{2N-M,M+1}=L_{2N-M,M+1}(0) 
$
with $\rank  L_{2N-M,M+1} = M$. That means that $L_{2N-M,M+1}$ has one non-zero kernel vector which is chosen to be the weight vector  $\boldsymbol{w}$ for the approximat (\ref{finbar}) (in case of exact data). In general case, $\boldsymbol{w}$ in (\ref{lls}) is computed as the right singular vector of the Loewner matrix $L_{2N+1-n,n}$ in (\ref{loew}) that corresponds to the smallest singular value of this matrix.

We refer the interested reader to \cite{DPP21} for details on how to solve the Loewner pencil (\ref{lp}). We note only that in case of exact data we use information regarding the numerical rank of $L_{2N-M,M+1}(t)$, $t=0,1$, and the eigenvalue problem (\ref{lp}) is reduced to a square one of the size $M\times M$.

\section{Rational structure of the Fourier coefficients of multivariate exponential sums}
\label{rfc}


We now address the problem of recovering of a multivariate exponential sum $f$ as in (\ref{mul}).
Analogously to the univariate case \cite{DP21}, we consider the  Fourier expansion of $f$  in $[0,P]^d$,
$$
f(\boldsymbol{t})= \sum \limits_{\boldsymbol{k}\in  \mathbb{Z}^d} c_{\boldsymbol{k}}(f) \mathrm{e}^{ \frac{2 \pi \mathrm{i}}{P}  \langle \boldsymbol{k}, \boldsymbol{t}\rangle}, \quad \boldsymbol{t} \in \rr^{d},
$$
with the Fourier coefficients $c_{\boldsymbol{k}}(f) \coloneqq \frac{1}{P^d}\int_{[0,P]^{d}} f (\boldsymbol{t}) \, \mathrm{e}^{-\frac{2 \pi \mathrm{i}}{P} \langle \boldsymbol{k},\boldsymbol{t}\rangle} \,  \mathrm{d} \boldsymbol{t}$,   $\boldsymbol{k} \in  \zz^d$.
Further, we will show that $c_{\boldsymbol{k}}(f)$ have a rational structure.   
Using (\ref{coef}),   we obtain
\begin{align}
c_{\boldsymbol{k}}(f) &  = \sum\limits_{j=1}^{M} \frac{\gamma_j}{P^d}   \int_{[0,P]^{d}}  \mathrm{e}^{\langle \boldsymbol{\lambda}_j, \boldsymbol{t} \rangle   - \frac{2 \pi \mathrm{i}}{P} \langle \boldsymbol{k},\boldsymbol{t}\rangle  }   \,  \mathrm{d} \boldsymbol{t} =  \sum\limits_{j=1}^{M}  \frac{\gamma_j}{P^d}  \int_{[0,P]^{d}}  \mathrm{e}^{\sum_{\ell=1}^{d}( \lambda_{j\ell} - 2 \pi \mathrm{i} k_\ell/P )t_\ell }  \,  \mathrm{d} \boldsymbol{t}  \notag\\
&=  \sum\limits_{j=1}^{M} \frac{\gamma_j}{P^d}  \prod\limits_{\ell=1}^{d}   \int\limits_0^{P}  \mathrm{e}^{(\lambda_{j\ell}-2 \pi \mathrm{i} k_\ell/P ) t_\ell} \,  \mathrm{d} t_\ell   =  \sum\limits_{j=1}^{M} \gamma_j  \prod\limits_{\ell=1}^{d}    \frac{\mathrm{e}^{\lambda_{j\ell}P}-1}{\lambda_{j\ell} P-2\pi \mathrm{i} k_\ell}  \label{coef2} 
\end{align}
for  $\lambda_{j\ell}\neq 2\pi \mathrm{i} k_\ell/P$. Note that in this paper we consider the case when all $\lambda_{j\ell}\neq 2\pi \mathrm{i} k_\ell/P$. If for some $j$ and $\ell$, $\lambda_{j\ell} =  2\pi \mathrm{i} k_\ell/P$, the formula (\ref{coef2})  will be modified according to (\ref{coef}) and the frequency components of this form can be reconstructed in a post-processing step (see Remark \ref{rem21}).
Using formula (\ref{coef2}),  we define a multivariate rational function 
\begin{equation}\label{rat}
r(\boldsymbol{z})=\sum\limits_{j=1}^{M}      \frac{a_{j}}{(z_1 - b_{j 1}) ...(z_d - b_{j d})}, \quad \boldsymbol{z}=(z_1,\ldots,z_d),
\end{equation}
 which, as follows from (\ref{coef2}), satisfies the interpolation conditions
 \begin{equation}\label{intd}
r(\boldsymbol{k})=c_{\boldsymbol{k}}(f)  ,  \quad \boldsymbol{k} \in \mathbb{Z}^{d},
\end{equation}
with parameters
\begin{align}
b_{j\ell}= \lambda_{j \ell} P/ 2\pi \mathrm{i} ,  \quad a_j= \frac{\gamma_j}{(2\pi \mathrm{i} )^d } \prod\limits_{\ell=1}^d (1- \mathrm{e}^{\lambda_{j \ell}P}).  \label{par}
\end{align}
That means that the multivariate exponential recovery problem (\ref{mul}) can be reformulated in terms of  the multivariate rational interpolation problem (\ref{intd}).

To reconstruct parameters of the multivariate rational function (\ref{rat}), we do not use  methods from the theory of multivariate rational functions (see Remark \ref{rmrf} for further explanations). Instead, inspired by the ideas of multivariate Prony-like methods developed in \cite{CL2018} and \cite{PT13mult}, we propose two approaches to reduce the multivariate rational interpolation problem to several univariate ones. These can then be solved using either the univariate AAA method for rational approximation, as described in Subsection \ref{ESPIRAwithAAA}, or the Loewner pencil approach presented in Subsection \ref{mpa}.





\section{Parameter estimation for bivariate exponential sums  }
\label{d2case}

We first consider the recovery of the exponential sum (\ref{mul}) for $d=2$, i.e. let a bivariate exponential sum of order $M$ be given by
\begin{equation}\label{mul2d}
f(\boldsymbol{t})=f(t_1,t_2)\coloneqq\sum\limits_{j=1}^{M} \gamma_j \mathrm{e}^{ \lambda_{j1}t_1+ \lambda_{j2}t_2 },
\end{equation}
with non-zero coefficients $\gamma_j \in \cc$ and frequency vectors $\boldsymbol{\lambda}_j=(\lambda_{j1},\lambda_{j2}) \in \cc^{2}$, $\boldsymbol{\lambda}_j \neq \boldsymbol{\lambda}_i$ for $i\neq j$. 
Our goal is to reconstruct the parameters $M \in \nn$, $\gamma_j \in \cc$ and $\boldsymbol{\lambda}_j=(\lambda_{j1},\lambda_{j2}) \in \cc^{2}$ for $j=1,...,M$, using as input information a finite set of the Fourier coefficients $c_{\boldsymbol{k}}(f)=c_{(k_1,k_2)}(f)$. Based on (\ref{coef2}), the coefficients  $c_{\boldsymbol{k}}(f)$ have the form
$$
c_{(k_1,k_2)}(f)   =  \sum\limits_{j=1}^{M} \gamma_j     \frac{(\mathrm{e}^{\lambda_{j1}P}-1)(\mathrm{e}^{\lambda_{j2}P}-1)}{(\lambda_{j1}P-2\pi \mathrm{i} k_1)(\lambda_{j2}P-2\pi \mathrm{i} k_2)}   
$$
and satisfy the  interpolation condition 
\begin{equation}\label{int2d}
r(k_1,k_2)=c_{(k_1,k_2)}(f)  ,  \quad \boldsymbol{k}=(k_1,k_2) \in \mathbb{Z}^{2},
\end{equation}
where $r(z_1,z_2)$ is a bivariate rational function  given by
\begin{equation}\label{ratf2d}
r(z_1,z_2)= \sum\limits_{j=1}^{M}     \frac{a_{j}}{(z_1 - b_{j1})(z_2 - b_{j2})},
\end{equation}
with $a_j= -\frac{\gamma_j}{4\pi^2}(1-\mathrm{e}^{\lambda_{j1}P}) (1-\mathrm{e}^{\lambda_{j2}P})$ and $b_{j\ell}= \lambda_{j \ell} P/ 2\pi \mathrm{i}$ for $\ell=1,2$.

\begin{algorithm}[ht]\caption{Sparse grid reconstruction of  bivariate exponential sums}
\label{alg1}
\small{
\textbf{Input:} $N \in \mathbb{N}$ large enough, $P>0$; \\
\phantom{\textbf{Input:}}  $\tau \in \mathbb{N}$, such that $ |\mathrm{Im}(\lambda_{jm})| < \frac{2\pi  \tau}{P}$, $m=1,2$, $j=1,...,M$; \\
\phantom{\textbf{Input:}} Fourier coefficients $c_{(k,0)}(f)$, $c_{(0,k)}(f)$, $k=-N,...,N$, and $c_{(k,k+2\tau)}(f)$, $k=-N,...,N-2\tau$;\\ 
\phantom{\textbf{Input:}} ${tol}>0$ tolerance for the approximation error within the AAA routine.
\begin{enumerate}
\item Compute $b_{j1}$ and $-a_j/b_{j2}$ as poles and coefficients, respectively,  of the univariate rational function $r$ in (\ref{int2d1}) using Fourier coefficients $c_{(k,0)}(f)=r(k,0)$, $k=-N,...,N$, and $b_{\sigma(j)2}$ as poles of the univariate rational function $r$ in (\ref{int2d11}) employing Fourier coefficients $c_{(0,k)}(f)=r(0,k)$, $k=-N,...,N$, via one of the methods from Section \ref{recuni}. The parameter $M$ is reconstructed as the number of iteration steps of AAA method minus 1;
\item Compute the vector $\boldsymbol{c}= (c_{11},...,c_{M1},c_{\sigma(1)2},...,c_{\sigma(M)2})^T$ as the least squares solution to the linear system  (\ref{sys}) using Fourier coefficients $c_{(k,k+2\tau)}(f)=r(k,k+2\tau)$, $k=-N,...,N-2\tau$. If for some $j$ and $k$, $c_{j 1}=- c_{k 2}$ and $-a_j/b_{j2}= c_{j1}+ c_{k2} \frac{b_{j1}}{b_{k2}}-2 \tau \frac{c_{j1}}{b_{k2}}$,   set $p(j)=k$.
\item Reconstruct the frequencies $\boldsymbol{\lambda}_j$ by 
$
\boldsymbol{\lambda}_j=(2\pi \mathrm{i} b_{j1}/P, 2\pi \mathrm{i} b_{\sigma(p(j))2}/P ), 
$ for $j=1,...,M$.
\item Compute $a_j$, $j=1,...,M$, by the formula $a_j=-(-a_j/b_{j2})b_{j2}$. Reconstruct coefficients $\gamma_j$, $j=1,...,M$ via (\ref{gamj}).
\end{enumerate}

\noindent
\textbf{Output:} $M$, $\gamma_{j}$, $\boldsymbol{\lambda}_j=(\lambda_{j1},\lambda_{j2})$, for $j=1, \ldots , M$ (all parameters of $f$).}
\end{algorithm}

\subsection{Sparse grid reconstruction of  bivariate exponential sums }
\label{2dcd}

In this subsection, we present our initial approach to recovering the bivariate exponential sum (\ref{mul2d}). Similarly to the univariate case, we begin by recovering the frequency vectors $\boldsymbol{\lambda}_j=(\lambda_{j1},\lambda_{j2})$ for $j=1,...,M$. This problem is equivalent to recovering the pole vectors $\boldsymbol{b}_j=(b_{j1}, b_{j2})$ for $j=1,...,M$ of the rational function (\ref{ratf2d}). We use the samples $r(k_1,k_2)=c_{(k_1,k_2)}(f)$ as input data and  present a method for selecting  indices $(k_1,k_2)$ of the Fourier coefficients from the 2D cube $[-N,N]^2\cap \zz^2$ with $N\geq M$, to reduce the bivariate rational interpolation problem (\ref{ratf2d}) to two univariate ones along the main coordinate directions.




For this subsection we have to assume that $b_{j m}\neq b_{i m}$ for $j \neq i$ and $m=1,2$. That means that all first (or second) components of a frequency  vector $\boldsymbol{\lambda}_j=(\lambda_{j1},\lambda_{j2})$ are pairwise distinct.  First, we  separately reconstruct two sets of poles $\{b_{j1} \}_{j=1}^M$ and $\{b_{j2} \}_{j=1}^M$. To achieve this, we use samples  of (\ref{ratf2d}) along the main coordinate directions, specifically at the points $r(k,0)$ and $r(0,k)$ for $k=-N,...,N$. Then we get two univariate rational interpolation problems, while the first one consists of the reconstruction of the poles $b_{j1}$ and the corresponding coefficients $-a_{j}/b_{j2}$ such that
\begin{equation}\label{int2d1}
r(k,0)= \sum\limits_{j=1}^{M}     \frac{-a_{j}/b_{j2}}{k - b_{j1}},  \quad k=-N,...,N,
\end{equation}
and the second one merely considers the computation of the poles $b_{j2}$ of the rational interpolant
\begin{equation}\label{int2d11}
r(0,k)= \sum\limits_{j=1}^{M}     \frac{-a_{j}/b_{j1}}{k - b_{j2}},  \quad k=-N,,...,N,
\end{equation}
particularly, the coefficients $-a_{j}/b_{j1}$, $j=1,...,M$ are not needed. 


It arises now the problematic that the order of the poles is not unique,
hence it can not be reconstructed by any reconstruction method. As a direct consequence we need to match the order of the poles of $r(k,0)$ and $r(0,k)$. Therefore, we assume that the poles $b_{j 1}$ of $r(k,0)$ are reconstructed in the correct order, while for $r(0,k)$ we have computed poles $b_{\sigma(j)2}$ for an unknown permutation $\sigma$. We now aim to find the correct pairing, that means the correct positions $p(j)$ of $j$ such that $\sigma(p(j))=j$.
In order to do it, we have to employ  additional samples of  (\ref{ratf2d}). Let $\tau \in \mathbb{N}$ be the smallest positive integer such that $|\mathrm{Im}(\lambda_{jm})| <   \frac{2\pi  \tau}{P}$  for $m=1,2$ and $j=1,...,M$ (that means  that  $ |\mathrm{Re} (b_{jm})| < \tau$) and we are also given the function values
\begin{equation}\label{ratkk1}
r(k,k+2\tau)= \sum\limits_{j=1}^{M}     \frac{a_{j}}{(k - b_{j1})(k - (b_{j2}-2\tau))}, \quad k=-N,...,N-2\tau,
\end{equation}
such that $-\tau < \mathrm{Re} (b_{j1}) < \tau$ and $-3 \tau <   \mathrm{Re} (b_{j2}) -2\tau < - \tau$. (We can always choose $N$ large enough such that $N-\tau\geq M$.)
Next, we compute the coefficients $c_{j1}$ and $c_{\sigma(j)2}$ of the partial fraction decomposition of (\ref{ratkk1}),
$$
r(k,k+2\tau)= \sum\limits_{j=1}^{M}     \frac{c_{j1}}{k - b_{j1}}+ \sum\limits_{j=1}^{M}     \frac{c_{\sigma(j)2}}{k - (b_{\sigma(j)2}-2\tau)}, \quad k=-N,...,N-2\tau.
$$
It can be done by finding the least squares solution of the following system of linear equations 
\begin{equation}\label{sys}
\left(\left.\mathbf{B}^{(1)}_{2N -2\tau +1, M} \right | \mathbf{B}^{(2)}_{2N -2\tau +1, M} \right) \, \boldsymbol{c}= \boldsymbol{r},
\end{equation}
with $\boldsymbol{c}= (c_{11},...,c_{M1},c_{\sigma(1)2},...,c_{\sigma(M)2})^T$, $\boldsymbol{r}=(r(-N,-N+2\tau),...,r(N-2\tau, N))^T$ and
$$
\mathbf{B}^{(1)}_{2N -2\tau +1, M}= \left( \frac{1}{k - b_{j1}} \right)_{k=-N,j=1}^{N-2\tau,M},  \quad \mathbf{B}^{(2)}_{2N -2\tau +1, M}= \left( \frac{1}{k - (b_{\sigma(j)2}-2\tau)} \right)_{k=-N,j=1}^{N-2\tau,M}.
$$
For  $\sigma(p(j))=j$ we have
\[
\resizebox{\textwidth}{!}{$
\begin{aligned}
 \sum\limits_{j=1}^{M}     \frac{c_{j1}}{k - b_{j1}}+ \sum\limits_{j=1}^{M}     \frac{c_{\sigma(p(j))2}}{k - (b_{\sigma(p(j))2}-2\tau)}  
      & = \sum\limits_{j=1}^{M} \frac{(c_{j1}+c_{\sigma(p(j))2}) k- (c_{j1} b_{\sigma(p(j))2}+c_{\sigma(p(j))2} b_{j1}-2\tau c_{j1})}{(k - b_{j1}) (k - (b_{\sigma(p(j))2}-2\tau))}, \notag \\
      & =\sum\limits_{j=1}^{M}     \frac{a_{j}}{(k - b_{j1})(k - (b_{\sigma(p(j))2}-2\tau))},
\end{aligned}
$}
\]
 which gives us the necessary conditions for the correct pairing
 \begin{align}
     c_{j1}=-c_{\sigma(p(j))2}, \quad
     a_j =-(c_{j1} b_{\sigma(p(j))2}+c_{\sigma(p(j))2} b_{j1}-2\tau c_{j1}).  \label{n2} 
 \end{align}
 Let us show that (\ref{n2}) are also the sufficient conditions. Let us assume that there is a different sequence $\rho(j)\neq p(j)$, i.e. $b_{\sigma(p(j))2} \neq b_{\sigma(\rho(j))2}$, such that $c_{j1}=-c_{\sigma(\rho(j))2}$ and $a_j =-(c_{j1} b_{\sigma(\rho(j))2}+c_{\sigma(\rho(j))2} b_{j1}-2\tau c_{j1})$. Then 
 $
c_{j1} b_{\sigma(p(j))2}+c_{\sigma(p(j))2} b_{j1}= c_{j1} b_{\sigma(\rho(j))2}+c_{\sigma(\rho(j))2} b_{j1}.
 $
 Since $c_{\sigma(p(j))2}=c_{\sigma(\rho(j))2}$, we get $c_{j1} b_{\sigma(p(j))2}=c_{j1} b_{\sigma(\rho(j))2}$. The vector $\boldsymbol{c}$ is a vector with nonzero coordinates (see Remark \ref{choicen}), thus we get the contradiction $ b_{\sigma(p(j))2}= b_{\sigma(\rho(j))2}$. To find the correct pairing, we use the coefficients $-a_{j}/b_{j2}$ computed together with poles $b_{j1}$ and we apply the following procedure: if for some $j$ and $k$, $c_{j 1}=- c_{k 2}$ and $-a_j/b_{j2}= c_{j1}+ c_{k2} \frac{b_{j1}}{b_{k2}}-2 \tau \frac{c_{j1}}{b_{k2}}$,  then set $p(j)=k$.  When we have the poles $\boldsymbol{b}_j=(b_{j1},b_{\sigma(p(j))2})$, we can determine the frequencies $\boldsymbol{\lambda}_j=\frac{2 \pi \mathrm{i}}{P}  \boldsymbol{b}_j$,  $j=1,...,M$.

We consider now the computation of coefficients $\gamma_j$. Using again the coefficients $-a_j/b_{j2}$, we compute $a_j$ by $a_j=-(-a_j/b_{j2})b_{j2}$, $j=1,...,M$ and then by 
\begin{equation}\label{gamj}
\gamma_j= - \frac{4\pi^2 a_j}{(1-\mathrm{e}^{\lambda_{j1}P}) (1-\mathrm{e}^{\lambda_{j2}P})}, \quad j=1,...,M,
\end{equation}
the wanted coefficients $\gamma_j$. 

The main ideas of this method are presented in Algorithm \ref{alg1}. We use in total $3(2N+1)-2\tau$ (i.e. $\mathcal{O}(N)$ ) Fourier coefficients. We employ $2(2N+1)$ Fourier coefficients $c_{(k,0)}(f)$ and  $c_{(0,k)}(f)$ for $k=-N,...,N$ to separately reconstruct the two sets of frequencies $\{\lambda_{j1} \}_{j=1}^M$ and $\{\lambda_{j2} \}_{j=1}^M$, and additional $2N+1-2\tau$ Fourier coefficients $c_{(k,k+2\tau)}(f)$ for $k=-N,...,N-2\tau$ to find the correct pairs $(\lambda_{j1},\lambda_{j2})$, $j=1,...,M$. Note that the corresponding indices $(k,0)$, $(0,k)$ and $(k,k+2\tau)$ are located on three straight lines in the 2D cube $[-N,N]^2$ (see Figure \ref{2dex}).
The complexity of the Algorithm \ref{alg1} is $\mathcal{O}(2 N M^3)$ and is determined by the complexity of the AAA method that we apply twice. Within the AAA routine, at each of $M+1$ iteration steps $n$ we compute  SVD of the Loewner matrix $L_{2N+1-n,n}$ for $n=1,...,M+1$ which requires  $\mathcal{O}(N M^2)$  flops.

\begin{remark}
\label{choicen}
We still have to explain the choice and the role of a parameter $\tau$. In case when $b_{j m}\neq b_{i \ell}$ for $j \neq i$ and $m \neq \ell$ we can choose $\tau=0$. But when we assume that we also have poles $b_{j m}= b_{i \ell}$ for $j \neq i$ and some $m$ and $\ell$, the choice $\tau=0$ will give us the matrix $(\mathbf{B}^{(1)}_{2N -2\tau +1, M}| \mathbf{B}^{(2)}_{2N -2\tau +1, M} )$ with at least two linearly dependent columns. To ensure that $(\mathbf{B}^{(1)}_{2N -2\tau +1, M}| \mathbf{B}^{(2)}_{2N -2\tau +1, M} )$ is of full rank and to guarantee the existence of a solution vector $\boldsymbol{c}$ with nonzero components,  we have to choose $\tau>0$, such that all $b_{j1}$ and $b_{j2}-2\tau$ are different. The choice of $\tau$ that satisfies  $|\mathrm{Re} (b_{jm})| < \tau$ solves the problem since then $-\tau < \mathrm{Re} (b_{j1}) < \tau$ and $-3\tau <   \mathrm{Re} (b_{j2}) -2\tau < - \tau$. For Algorithm \ref{alg1} we have to assume that we know the segment $[-\frac{2\pi \tau}{P},\frac{2\pi \tau}{P}]$ at which the imaginary parts of the frequency components are located.

\end{remark}




\begin{algorithm}[ht]\caption{Recovery of bivariate exponential sums via recursive dimension reduction }
\label{alg2}
\small{
\textbf{Input:} $N \in \mathbb{N}$ large enough, $P>0$; \\
\phantom{\textbf{Input:}} Fourier coefficients $c_{(k,i)}(f)$,  $k,i=-N,...,N$;\\ 
\phantom{\textbf{Input:}} ${tol}>0$ tolerance for the approximation error within the AAA routine.

\begin{enumerate}
\item Compute  $b_{j1}$, $j=1,...,M_1$   as pairwise distinct poles  of the univariate rational function $r$ in (\ref{unr}) via one of the methods from Section \ref{recuni} using the samples $c_{(k,0)}(f)=r(k,0)$, $k=-N,...,N$. The number $M_1$ is reconstructed as number of iteration steps of the AAA method minus $1$.
\item For each $m=1,...,M_1$, compute  $d_m(i)$, $i=-N,...,N$, as the least squares solutions to the   linear systems (\ref{dl}) using Fourier coefficients $c_{(k,i)}(f)=r(k,i)$, $k=-N,...,N$.
\item For each $m=1,...,M_1$, using the values $d_m(i)$, $i=-N,...,N$, compute the poles $b_{\ell 2}^{(j_m)}$, $\ell=1,...,\mu_m$ of a rational function 
(\ref{dl11}) applying again one of the method from Section \ref{recuni}. For each $m=1,...,M_1$,  reconstruct parameters $\mu_m$ as numbers of iteration steps of the AAA algorithms minus 1. Finally, reconstruct $M$ as $\sum_{m=1}^{M_1} \mu_m=M$.
\item Reconstruct the frequency vectors by
$
\boldsymbol{\lambda}_{m, \ell}=(2\pi \mathrm{i} b_{j_m1}/P, 2\pi \mathrm{i} b_{\ell 2}^{(j_m)} /P)$, for $\ell=1,...,\mu_m$, $m=1,...,M_1$.

\item Compute $a_j$, $j=1,...,M$, as the least squares solution to the linear system (\ref{sys2}) using Fourier coefficients $c_{(k,i)}(f)$, $k,i=-N,...,N$. Reconstruct coefficients $\gamma_j$ via (\ref{gamj}).
\end{enumerate}

\noindent
\textbf{Output:} $M$, $\gamma_{j}$, $\boldsymbol{\lambda}_j=(\lambda_{j1},\lambda_{j2})$, for $j=1, \ldots , M$ (all parameters of $f$).}
\end{algorithm}

\subsection{Recovery of bivariate exponential sums via recursive dimension reduction }
\label{rec2d}
In this subsection, we describe the second method for the recovery of the bivariate exponential sum (\ref{mul2d}). Note that this method works without the assumption $\lambda_{j m} \neq \lambda_{i m}$ for $j \neq i$ and $m=1,2$, which was required  in Subsection \ref{2dcd}.

Let us assume that the rational function (\ref{ratf2d}) has $M_1\leq M$ pairwise distinct poles $b_{j_m 1}$ with  multiplicities  $\mu_m$ for $m=1,...,M_1$ (i.e. $\mu_m$ is the number of poles $b_{j 1}$ equal to   $b_{j_m 1}$). Then  $\sum_{m=1}^{M_1} \mu_m=M$ and the rational function (\ref{ratf2d}) can be rewritten as
$$
r(z_1,z_2)=\sum\limits_{m=1}^{M_1}  \frac{1}{z_1 - b_{j_m 1}}  \sum\limits_{\ell=1}^{\mu_m}  \frac{a_\ell^{(j_m)}}{z_2-b_{\ell  2}^{(j_m)}} = \sum\limits_{m=1}^{M_1}  \frac{d_m(z_2)}{z_1 - b_{j_m 1}} ,
$$
where the poles $b_{\ell 2}^{(j_m)}$, $\ell=1,...,\mu_m$ are pairs to the pole $b_{j_m 1}$,   $a_{\ell }^{(j_m)}$ is the weight related to the pole $b_{\ell  2}^{(j_m)}$, and for each $m$, $d_m(z_2)=\sum_{\ell=1}^{\mu_m}  \frac{a_\ell^{(j_m)}}{z_2-b_{\ell  2}^{(j_m)}}$ is a rational function itself of type $(\mu_m-1,\mu_m)$ with poles $b_{\ell 2}^{(j_m)}$.

Applying the AAA routine from Subsection \ref{ESPIRAwithAAA} or the Loewner pencil from Subsection \ref{mpa} with
 samples 
\begin{equation}\label{unr}
r(k,0)=\sum\limits_{m=1}^{M_1}  \frac{d_m(0)}{k - b_{j_m1}} \quad  \text{ where } \quad d_m(0)= -\sum\limits_{\ell=1}^{\mu_m}  \frac{a_\ell^{(j_m)}}{b_{\ell 2}^{(j_m)}},
\end{equation}
for $k=-N,...,N$, we reconstruct pairwise distinct poles $b_{j_m1}$, $m=1,...,M_1$, and the parameter $M_1$.
We still want to reconstruct  $M$ and $b_{\ell 2}^{(j_m)}$ for each $m$.  We notice that $\mu_m\leq M-M_1 +1 \leq N-M_1+1 \leq N$. Further, for  each $i=-N,...,N$, we compute the coefficients $d_m(i)$ as the  least square solution of the systems of linear equations
\begin{equation}\label{dl}
    r(k,i)=\sum\limits_{m=1}^{M_1}  \frac{d_m(i)}{k - b_{j_m1}}, \quad k=-N,...,N.
\end{equation}
For each $m=1,...,M_1$, having values $d_m(i)$ of  $d_m$ at  $i=-N,..., N$, 
\begin{equation}\label{dl11}
d_m(i)=\sum\limits_{\ell=1}^{\mu_m}  \frac{a_\ell^{(j_m)}}{i-b_{\ell 2}^{(j_m)}}
\end{equation}
we can reconstruct $\mu_m$ poles $b_{\ell 2}^{(j_m)}$, $\ell=1,...,\mu_m$, by applying one of the methods from Section \ref{recuni}. Note that using this approach, we automatically solve the problem with the correct pairing, since in this case we reconstruct the poles $b_{\ell 2}^{(j_m)}$, $\ell=1,...,\mu_m$ already as  correct pairs to the pole $b_{j_m 1}$ for each $m=1,...,M_1$.
Finally, we reconstruct the frequency vectors $\boldsymbol{\lambda}_j$ as 
$
\boldsymbol{\lambda}_{m,\ell}=(2\pi \mathrm{i} b_{j_m1}/P, 2\pi \mathrm{i} b_{\ell 2}^{(j_m)} /P)$, $\ell=1,...,\mu_m$, $m=1,...,M_1$,
where $\sum_{m=1}^{M_1} \mu_m=M$. 




Then we compute coefficients $a_j$, $j=1,...,M$, as the least squares solution of the system of linear equations 
\begin{equation}\label{sys2}
r(k, \ell)= \sum\limits_{j=1}^{M}     \frac{a_{j}}{(k - b_{j1})( \ell - b_{j2})}, \quad k,\ell=-N,...,N, 
\end{equation}
and  reconstruct coefficients $\gamma_j$, $j=1,...,M$, by  (\ref{gamj}). 

The main ideas of this method are described in Algorithm \ref{alg2}. As input data we apply  $(2N+1)^2$ (i.e. $\mathcal{O}(N^2)$ ) Fourier coefficients  $c_{\boldsymbol{k}}(f)$, $ \boldsymbol{k}\in [-N,N]^2\cap \zz^2$.  Algorithm \ref{alg2} has an overall computational cost of $\mathcal{O}(N M^2(M^2+N))$, which is determined by the AAA method, which we apply $1+M_1\leq 1+M$ times, and  the solutions to the linear systems (\ref{dl}). First, we reconstruct the pairwise distinct frequencies $\lambda_{j_m 1}$ applying AAA, which requires $\mathcal{O}(N M^3)$ flops (we compute SVD of the Loewner matrices $L_{2N+1-n,n}$ for $n=1,...,M_1+1$). Then we apply the AAA method $M_1\leq M$ times to reconstruct the frequencies $\lambda_{\ell 2}^{(j_m)}$, $\ell=1,...,\mu_m, \, m=1,...,M_1$,
where $\sum_{m=1}^{M_1} \mu_m=M$.  Here for each $m=1,...,M_1$ we compute SVD of Loewner matrices $L_{2N+1-n,n}$  for $n=1,..,\mu_m+1$, which requires  $\mathcal{O}(N M^2)$ flops. The computational cost of the solution of linear systems (\ref{dl}) is  $\mathcal{O}(N^2 M^2)$.


\begin{remark}\label{rmrf}
The two  recovery approaches for bivariate exponential sums described above are based on the idea of reducing a bivariate rational interpolation problem to several univariate ones. Of course, the direct idea will be to reconstruct bivariate exponential sums by applying bivariate rational interpolation techniques in the Fourier domain. Using access  to the Fourier coefficients (some finite set) given by the interpolation property (\ref{int2d}), we can compute a bivariate rational function in the barycentric representation
\begin{equation}\label{2dbar}
R(z_1,z_2) = \sum\limits_{s=1}^{M+1}  \sum\limits_{i=1}^{M+1} \frac{c_{(k_{s},p_i)}(f)w_{s i}}{(z_1-k_s)(z_2-p_i) } \bigg/ \sum\limits_{s=1}^{M+1}  \sum\limits_{i=1}^{M+1} \frac{w_{s i}}{(z_1-k_s)(z_2-p_i) } 
\end{equation}
via the $p$-AAA algorithm developed in \cite{R22}. To reconstruct the poles $(b_{j1},b_{j2})$ we have to compute the zeros of the denominator of (\ref{2dbar}). To our knowledge, there is no  analog of the generalized eigenvalue problem (\ref{eig}) in the bivariate case.
 
\end{remark}

\section{Recovery of $d$-variate exponential sums for $d>2$}
\label{spd}

Let us now consider the question of recovery of a $d$-variate, $d>2$,  exponential sums (\ref{mul}), using as input data its Fourier coefficients  
\begin{equation}\label{fcd}
   c_{\boldsymbol{k}}(f) =  \sum\limits_{j=1}^{M} \frac{a_j}{(k_1-b_{j1})\ldots (k_d-b_{j d})}, \quad \boldsymbol{k}=(k_1,...,k_d) \in \zz^d,
\end{equation}
with $b_{j\ell}$ and $a_j$ as in (\ref{par}). Recall that a multivariate rational function $r$ as in (\ref{rat}) with parameters (\ref{par}) satisfies the interpolation property (\ref{intd}) with Fourier coefficients (\ref{fcd}).


\begin{algorithm}[t]\caption{Sparse grid reconstruction of  $d$-variate exponential sums}
\label{alg3}
\small{
\textbf{Input:} $N \in \mathbb{N}$ large enough, $P>0$; \\
\phantom{\textbf{Input:}}  $\tau \in \mathbb{N}$, such that $|\mathrm{Im}(\lambda_{j m})| < \frac{2\pi  \tau}{P}$, $j=1,...,M$, $m=1,...,d$; \\
\phantom{\textbf{Input:}} Fourier coefficients $c_{(\boldsymbol{0}_{m-1},k,\boldsymbol{0}_{d-m})}(f)$ for $m=1,...,d$ and $k=-N,...,N$, \\ 
\phantom{\textbf{Input:}}  and   $c_{(\boldsymbol{0}_{m-2},k,k+2\tau,\boldsymbol{0}_{d-m})}(f)$  for $m=2,...,d$ and $k=-N,...,N-2\tau$;\\ 
\phantom{\textbf{Input:}} ${tol}>0$ tolerance for the approximation error within the AAA routine.

\begin{enumerate}
\item Compute the sets of poles  $\{b_{j1} \}_{j=1}^M$, $\{b_{\sigma_{m-1}(j)m} \}_{j=1}^M$, $2 \leq m \leq d$ and coefficients $A_{j m}$, $j=1,...,M$, $1\leq m \leq d$, using Fourier coefficients $c_{(\boldsymbol{0}_{m-1},k,\boldsymbol{0}_{d-m})}(f)=r(\boldsymbol{0}_{m-1},k,\boldsymbol{0}_{d-m})$, $k=-N,...,N$, as  
 poles and coefficients, respectively, of the univariate rational function (\ref{unratd})  via methods described in Section \ref{recuni}.
\item To find the correct pairing for  $\{b_{j1} \}_{j=1}^M$, $\{b_{\sigma_{m-1}(j)m} \}_{j=1}^M$, $2 \leq m \leq d$ use the following approach: for each $m=2,...,d$, assuming that $\{b_{\sigma_{m-2}(p_{m-2}(j))(m-1)} \}_{j=1}^M$ has the correct order, i.e. the sequence $p_{m-2}(j)$  is known already ($p_0\equiv id$, $\sigma_0\equiv id$), find the correct order for $\{b_{\sigma_{m-1}(j)m} \}_{j=1}^M$, i.e. the sequence $p_{m-1}(j)$ such that $\sigma_{m-1}(p_{m-1}(j))=j$. Compute the vector   $\boldsymbol{c}^{(m-1)}= (c_{11}^{(m-1)},...,c_{M1}^{(m-1)},c_{\sigma_{m-1}(1)2}^{(m-1)},...,c_{\sigma_{m-1}(M)2}^{(m-1)})^T$ as the least squares solution to the linear system  (\ref{sysd}), using Fourier coefficients $c_{(\boldsymbol{0}_{m-2},k,k+2\tau,\boldsymbol{0}_{d-m})}(f)=r(\boldsymbol{0}_{m-2},k,k+2\tau,\boldsymbol{0}_{d-m})$, $k=-N,...,N-2\tau$.
 If for some $j$ and $k$, $   c_{j1}^{(m-1)} =- c_{k2}^{(m-1)}$ and 
 $A_{j (m-1)}= c_{j1}^{(m-1)}+ c_{k2}^{(m-1)} \frac{b_{j(m-1)}}{b_{k m}}-2 \tau \frac{c_{j1}^{(m-1)}}{b_{k m}}$,  set $p_{(m-1)}(j)=k$

\item Reconstruct frequencies  
$
\boldsymbol{\lambda}_j=(2\pi \mathrm{i} b_{j1}/P, 2\pi \mathrm{i} b_{\sigma_1(p_1(j))2}/P,..., 2\pi \mathrm{i} b_{\sigma_{d-1}(p_{d-1}(j))d}/P)$, $j=1,...,M.
$
\item Compute $a_j$, $j=1,...,M$ by  $a_{j}= A_{j 1}  \prod_{\ell=2 }^d (-b_{j \ell})$. Reconstruct coefficients $\gamma_j$ by  (\ref{gamjd}).
\end{enumerate}

\noindent
\textbf{Output:} $M$, $\gamma_{j}$, $\boldsymbol{\lambda}_j=(\lambda_{j1},...,\lambda_{j d})$, for $j=1, \ldots , M$ (all parameters of $f$).}
\end{algorithm}

\subsection{Sparse grid reconstruction of  $d$-variate exponential sums }

Similarly as for $d=2$, we start with the reconstruction of the frequency vectors $\boldsymbol{\lambda}_j=(\lambda_{j1},\ldots, \lambda_{jd})$, $j=1,...,M$. That means that first we have to reconstruct the pole vectors $\boldsymbol{b}_j=(b_{j1},...,b_{j d})$ of the rational function (\ref{rat}) using as input data $c_{\boldsymbol{k}}(f)$ for $\boldsymbol{k} \in [-N,N]^d \cap  \zz^d$. Afterwards,  the relation 
$
    \boldsymbol{\lambda}_j=(2\pi \mathrm{i} b_{j1}/P,..., 2\pi \mathrm{i} b_{j d}/P )$, $j=1,..,M$,
will give us the wanted frequency vectors.

We adapt the 2-variate ideas from Subsection \ref{2dcd} to the case of arbitrary dimension $d \in \nn$. Further, we use the notation  $ \boldsymbol{0}_{\ell}=(0,...,0)$ to denote the zero vector of the length $\ell$.
First, we  separately reconstruct sets of poles $\{b_{j1} \}_{j=1}^M$, ..., $\{b_{j d} \}_{j=1}^M$ and then we find the correct pairing. We work under the assumption that $b_{j m}\neq b_{i m}$ for $j \neq i$ and $m=1,...,d$. To compute the sets  $\{b_{j m} \}_{j=1}^M$ for $1\leq m \leq d$, we use samples of $r(z)$ as in (\ref{rat}) along all $d$ coordinate directions, i.e. for $k=-N,...,N$ we consider
\begin{equation}\label{unratd}
r(\boldsymbol{0}_{m-1},k,\boldsymbol{0}_{d-m})=\sum\limits_{j=1}^{M} \frac{A_{j m}}{k-b_{j m}},\quad \text{ with }\quad A_{j m}= a_{j} \bigg/ \prod\limits_{\substack{\ell=1 \\  \ell\neq m } }^d (-b_{j \ell}) .
\end{equation}
 Together with the poles $\{b_{j m} \}_{j=1}^M$ we compute also the coefficients $A_{j m}$, $j=1,...,M$ (coefficients $A_{j d}$ are not needed). 

Further, we explain how we find the correct pairing.  We assume  that the poles $\{b_{j1}\}_{j=1}^M$ have the correct order and the poles $\{b_{\sigma_{m-1}(j)m} \}_{j=1}^M$, $2 \leq m \leq d$ are permuted and we have to find the sequences $p_{m-1}(j)$ such that  $\sigma_{m-1}(p_{m-1}(j))=j$. We apply the following idea: assuming that the poles $\{b_{j1}\}_{j=1}^M$ have the correct order we find the sequence $p_1(j)$ and the correct   pairing $(b_{j1}, b_{\sigma_1(p_1(j))2})$, then  the poles $\{b_{\sigma_1(p_1(j))2})\}_{j=1}^M$ have the correct order and we find $p_2(j)$ which yields the correct   pairing $(b_{\sigma_1(p_1(j))2}), b_{\sigma_2(p_2(j))3})$, etc. and at the end when the  poles $\{b_{\sigma_{d-2}(p_{d-2}(j))(d-1)}\}_{j=1}^M$ have the correct order we find $p_{d-1}(j)$ and the correct   pairing $(b_{\sigma_{d-2}(p_{d-2}(j))(d-1)}, b_{\sigma_{d-1}(p_{d-1}(j))d})$. Thus, the wanted pole vectors $\boldsymbol{b}_j$ are reconstructed as $\boldsymbol{b}_j=(b_{j1}, b_{\sigma_1(p_1(j))2},..., b_{\sigma_{d-1}(p_{d-1}(j))d})$, $j=1,...,M$.

Let again $\tau \in \mathbb{N}$ be the smallest positive integer such that $|\mathrm{Im}(\lambda_{j m})| < \frac{2\pi  \tau}{P}$  for all $j=1,...,M$ and $m=1,...,d$ (it means that also $|\mathrm{Re} (b_{j m})| < \tau$).
To find the correct pairing for $\{b_{\sigma_{m-2}(p_{m-2}(j))(m-1)} \}_{j=1}^M$ and  $\{b_{\sigma_{m-1}(j)m} \}_{j=1}^M$, $2 \leq m \leq d$,   assuming that the poles $\{b_{\sigma_{m-2}(p_{m-2}(j))(m-1)} \}_{j=1}^M$ (we denote $p_0\equiv id$, $\sigma_0=id$) have  correct order and we have to find the sequence $p_{m-1}(j)$, we use the following approach. For simplicity, we use that already  $\sigma_{m-2}(p_{m-2}(j))=j$. We employ the samples 
\begin{equation}\label{pdf3}
\resizebox{.92\textwidth}{!}{%
$r(\boldsymbol{0}_{m-2},k,k+2\tau, \boldsymbol{0}_{d-m}) = \sum\limits_{j=1}^{M} \frac{\hat{a}_j^{(m-1)}}{(k-b_{j(m-1)}) (k-(b_{j m}-2\tau))}, \quad \text{with}\quad \hat{a}_j^{(m-1)}= \frac{a_{j}} { \prod\limits_{\substack{\ell=1 \\  \ell\neq m-1,m } }^d (-b_{j \ell}) } $}
\end{equation}
for $k=-N,...,N-2\tau$. According to our assumption  $ |\mathrm{Re} (b_{j(m-1)})| < \tau$ as well as $-3\tau < \mathrm{Re} (b_{j m}) - 2\tau< -\tau$.
Next, we compute the vector   $\boldsymbol{c}^{(m-1)}= (c_{11}^{(m-1)},...,c_{M1}^{(m-1)},$ $c_{\sigma_{m-1}(1)2}^{(m-1)},...,c_{\sigma_{m-1}(M)2}^{(m-1)})^T$ 
of the partial fraction decomposition of (\ref{pdf3}),
$$
r(\boldsymbol{0}_{m-2},k,k+2\tau, \boldsymbol{0}_{d-m}) = \sum\limits_{j=1}^{M}     \frac{c_{j1}^{(m-1)}}{k - b_{j(m-1)}}+ \sum\limits_{j=1}^{M}     \frac{c^{(m-1)}_{\sigma_{m-1}(j)2}}{k - (b_{\sigma_{m-1}(j)m}-2\tau)}, 
$$
$k=-N,...,N-2\tau$, as the least squares solution of the linear system
\begin{equation}\label{sysd}
(\mathbf{B}_{2N-2\tau+1,M}^{(m-1)}| \mathbf{B}_{2N-2\tau+1,M}^{(m-1)} ) \boldsymbol{c}^{(m-1)}= \boldsymbol{r}^{(m-1)},
\end{equation}
 with 
\[\resizebox{\textwidth}{!}{%
$\mathbf{B}_{2N-2\tau+1,M}^{(m-1)}= \left( \frac{1}{k - b_{j(m-1)}} \right)_{k=-N,j=1}^{N-2\tau,M},  \quad \mathbf{B}_{2N-2\tau+1,M}^{(m-1)}= \left( \frac{1}{k - (b_{\sigma_{m-1}(j)m}-2\tau)} \right)_{k=-N,j=1}^{N-2\tau,M}
$}\]
and 
$
\boldsymbol{r}^{(m-1)}=r(\boldsymbol{0}_{m-1},-N,-N+2\tau, \boldsymbol{0}_{d-m-1}), ...,r(\boldsymbol{0}_{m-1},N-2\tau, N, \boldsymbol{0}_{d-m-1}))^T.
$
Since  for  $\sigma_{m-1}(p_{m-1}(j))=j$ we have
\begin{align*}
  \sum\limits_{j=1}^{M}  \frac{\hat{a}_j^{(m-1)}}{(k-b_{j(m-1)}) (k-(b_{j m}-2\tau))} &=\sum\limits_{j=1}^{M}     \frac{c_{j1}^{(m-1)}}{k - b_{j(m-1)}}+ \sum\limits_{j=1}^{M}     \frac{c^{(m-1)}_{\sigma_{m-1}(p_{m-1}(j))2}}{k - (b_{\sigma_{m-1}(p_{m-1}(j))m}-2\tau)}   \\
  & =\sum\limits_{j=1}^{M} \frac{B_{j(m-1)}}{(k - b_{j(m-1)}) (k - (b_{\sigma_{m-1}(p_{m-1}(j))m}-2\tau))},
\end{align*}
with
\begin{align*}
  B_{j(m-1)} &=(c_{j1}^{(m-1)}+c_{\sigma_{m-1}(p_{m-1}(j))2}^{(m-1)}) k \\
  & - (c_{j1}^{(m-1)} b_{\sigma_{m-1}(p_{m-1}(j))m}+c_{\sigma_{m-1}(p_{m-1}(j))2} ^{(m-1)} b_{j(m-1)}-2\tau c_{j1}^{(m-1)}) , 
\end{align*}
and applying similar ideas as in Subsection \ref{2dcd}, we get the necessary  and sufficient conditions for the correct pairing:
\begin{align*}
    c_{j1}^{(m-1)} &=- c_{\sigma_{m-1}(p_{m-1}(j))2}^{(m-1)}, \\
    \hat{a}_j^{(m-1)} & = - (c_{j1}^{(m-1)} b_{\sigma_{m-1}(p_{m-1}(j))m}+c_{\sigma_{m-1}(p_{m-1}(j))2} ^{(m-1)} b_{j(m-1)}-2\tau c_{j1}^{(m-1)}).
\end{align*}
Recall that
$$A_{j (m-1)}= a_{j} \bigg/ \prod\limits_{\substack{\ell=1 \\  \ell\neq m-1 } }^d (-b_{j \ell}) = -\frac{\hat{a}_j^{(m-1)}}{b_{j m}}.$$ If for some $j$ and $k$, $   c_{j1}^{(m-1)} =- c_{k2}^{(m-1)}$ and 
 $A_{j (m-1)}= c_{j1}^{(m-1)}+ c_{k2}^{(m-1)} \frac{b_{j(m-1)}}{b_{k m}}-2 \tau \frac{c_{j1}^{(m-1)}}{b_{k m}}$, then set $p_{(m-1)}(j)=k$. Finally, we reconstruct the frequencies $\boldsymbol{\lambda}_j$ for $j=1,...,M$ by
$
\boldsymbol{\lambda}_j=(2\pi \mathrm{i} b_{j1}/P, 2\pi \mathrm{i} b_{\sigma_1(p_1(j))2}/P,..., 2\pi \mathrm{i}b_{\sigma_{d-1}(p_{d-1}(j))d}/P)$.


We consider now the computation of the coefficients $\gamma_j$. First, we reconstruct the coefficients $a_j$, $j=1,...,M$ by 
$a_{j}= A_{j 1}  \prod_{\ell=2 }^d (-b_{j \ell}) $
 and then the wanted coefficients $\gamma_j$ by
\begin{equation}\label{gamjd}
\gamma_j=  \frac{ a_j (2\pi \mathrm{i} )^d}{\prod\limits_{\ell=1}^d (1- \mathrm{e}^{\lambda_{j \ell} P})}, \quad j=1,...,M.
\end{equation}

 In Algorithm \ref{alg3}, we describe the main ideas of this method.
Similarly to the 2D case, we do not use here the full grid $ \boldsymbol{k} \in [-N,N]^d \cap \zz^d$  of the Fourier coefficients $c_{\boldsymbol{k}}(f)$. Instead, we chose a sparse grid of indices   $\boldsymbol{k}$  such that we reduce the $d$-variate rational approximation problem (\ref{fcd}) to $d$ univariate ones (\ref{unratd}), which are solved by a univariate method for rational approximation described in Section \ref{recuni}.  First,  we employ $d(2N+1)$  Fourier coefficients  $c_{(\boldsymbol{0}_{m-1},k,\boldsymbol{0}_{d-m})}(f)$ for $k=-N,...,N$ and $m=1,...,d$, to recover the frequency components $\{\lambda_{j1}\}$,...,$\{\lambda_{j d}\}$ for $j=1,...,M$. Further, we employ additional $(d-1)(2N+1-2\tau)$  Fourier coefficients $c_{(\boldsymbol{0}_{m-2},k,k+2 \tau,\boldsymbol{0}_{d-m})}(f)$ for $k=-N,...,N-2\tau$ and $m=2,...,d$, in order to create the correct pairs $(\lambda_{j1},\ldots, \lambda_{j d})$, $j=1,...,M$. For example, for $d=3$  the corresponding indices $(k,0,0)$, $(0,k,0)$, $(0,0,k)$, $(k,k+2\tau,0)$ and $(0,k,k+2\tau)$ are located on five straight lines in the 3D cube $[-N,N]^3$ (see Figure \ref{2dex1}).
In general, the method uses $\mathcal{O}(d N)$ Fourier coefficients and has the overall computational cost of $\mathcal{O}(d N M^3)$, which is determined by AAA method that we apply $d$ times.

\begin{algorithm}[h]\caption{Recovery of $d$-variate  exponential sums via recursive dimension reduction }
\label{alg4}
\small{
\textbf{Input:} $N \in \mathbb{N}$ large enough, $P>0$; \\
\phantom{\textbf{Input:}} Fourier coefficients $c_{(k_1,...,k_d)}(f)$,   $(k_1,...,k_d) \in [-N,N]^{d} \cap \zz^{d}$;\\ 
\phantom{\textbf{Input:}} ${tol}>0$ tolerance for the approximation error within the AAA routine.

\begin{enumerate}
\item Compute  $b_{j_m1}$, $m=1,...,M_1$   as pairwise distinct poles  of  $r$ in (\ref{rat4}) via one of the methods  from Section \ref{recuni} employing samples $c_{(k_1,\boldsymbol{0}_{d-1})}(f)=r(k_1,\boldsymbol{0}_{d-1})$, $k_1=-N,...,N$. The number $M_1$ is reconstructed as number of iteration steps of the AAA method minus $1$.
\item 
For each $m=1,...,M_1$, compute values of $d_{m}$ as in  (\ref{dminus1r})  on the grid of  points $(k_2,...,k_d) \in [-N,N]^{d-1} \cap \zz^{d-1}$ as least squares solutions of  linear systems (\ref{dl1}) using Fourier coefficients $c_{(k_1,...,k_d)}(f)=r(k_1,...,k_d)$,  $(k_1,...,k_d) \in [-N,N]^{d} \cap \zz^{d}$.
\item Compute pairwise distinct poles $b_{\ell_{s_1}^{(1)} 2}^{(j_m)}$, $s_1=1,...,M_2^{(j_m)}$, for each $m=1,...,M_1$, of a univariate rational function (\ref{dm1}) applying one of the methods  from Section \ref{recuni} using samples $d_{m}(k_2,\boldsymbol{0}_{d-2})$, $k_2=-N,...,N$.
 Compute  values of $d_{s_1}^{(m)}$ as in (\ref{ds1}) on the grid $(k_3,...,k_d) \in [-N,N]^{d-2} \cap \zz^{d-2}$ as solutions to the linear systems (\ref{dm1s}) using values $d_{m}(k_2, k_3,...,k_d)$ on the grid $(k_2,...,k_d) \in [-N,N]^{d-1} \cap \zz^{d-1}$.
\item For each step $p=3,...,d$, compute pairwise distinct poles  $b_{\ell_{s_{p-1}}^{(p-1)},p}^{(j_m,\ell_{s_1}^{(1)},...,\ell_{s_{p-2}}^{(p-2)})}$, $s_{p-1}=1,...,M_p^{(j_m,\ell_{s_1}^{(1)},...,\ell_{s_{p-2}}^{(p-2)})}$ of  univariate rational functions (\ref{dsp2}) applying one of the methods  from Section \ref{recuni} using the function values $d_{s_{p-2}}^{(m,s_1,...,s_{p-3})}(k_p,\boldsymbol{0}_{d-p})$, $k_p=-N,...,N$. If $p\leq d-1$ compute  values of $d_{s_{p-1}}^{(m,s_1,...,s_{p-2})}$ as in (\ref{dp1r}) on the grid $(k_{p+1},...,k_d) \in [-N,N]^{d-p} \cap \zz^{d-p}$ as solution of systems of linear equations (\ref{dsp2s}) using values $ d_{s_{p-2}}^{(m,s_1,...,s_{p-3})}(k_p,...,k_d)$ on the grid $(k_{p},...,k_d) \in [-N,N]^{d-p+1} \cap \zz^{d-p+1}$.
\item
Reconstruct poles $\boldsymbol{b}_{m,s_1,...,s_{d-1}}$ via (\ref{pol1}) and the  frequencies $\boldsymbol{\lambda}_j$ as vectors $\boldsymbol{\lambda}_{m,s_1,...,s_{d-1}}=\frac{2\pi \mathrm{i}}{P} \boldsymbol{b}_{m,s_1,...,s_{d-1}}$
for $m=1,...,M_1$, $s_1=1,...,M_2^{(j_m)}$, $s_{p-1}=1,...,M_p^{(j_m,\ell_{s_1}^{(1)},...,\ell_{s_{p-2}}^{(p-2)})}$, $p=3,...,d$. The parameter $M$ is reconstructed via (\ref{mrec}) with $p=d$.
\item Compute $a_j$, $j=1,...,M$, as a least squares solution to  the linear system (\ref{acoef}) using Fourier coefficients $c_{\boldsymbol{k}}(f)=r(\boldsymbol{k})$, $\boldsymbol{k} \in [-N,N]^d$. Reconstruct  $\gamma_j$ by  (\ref{gamjd}).
\end{enumerate}

\noindent
\textbf{Output:} $M$, $\gamma_{j}$, $\boldsymbol{\lambda}_j=(\lambda_{j1},...,\lambda_{j d})$, for $j=1, \ldots , M$ (all parameters of $f$).}
\end{algorithm}
\subsection{Recovery of $d$-variate exponential sums via  recursive dimension reduction }

We now generalize the bivariate method from Subsection \ref{rec2d} to any dimension $d> 2$. Let us have function values $r(k_1,...,k_d)$ of the $d$-variate rational function (\ref{rat}) on the grid $(k_1,...,k_d) \in [-N,N]^{d} \cap \zz^{d}$.  We rewrite this function as follows
$$
    r(\boldsymbol{z})= \sum\limits_{m=1}^{M_1}  \frac{1}{(z_1 - b_{j_m 1})} \sum\limits_{\ell=1}^{\mu_{m}}    \frac{a_{\ell}^{(j_m)}}{(z_2 - b_{\ell 2}^{(j_m)}) ...(z_d - b_{\ell d}^{(j_m)})} =  \sum\limits_{m=1}^{M_1}  \frac{d_{m}(z_2,...,z_d)}{z_1 - b_{j_m 1}},
    $$
where all poles $b_{j_m 1}$, $m=1,...,M_1$ are pairwise distinct, $\mu_{m}$ describes the multiplicity of $b_{j_m 1}$, i.e. $\sum_{m=1}^{M_1} \mu_m=M$, $a_{\ell}^{(j_m)}$ are some coefficients, and for each $m=1,...,M_1$,  $d_{m}(z_2,...,z_d)$ is a $(d-1)$-variate rational function of the form
\begin{equation}\label{dminus1r}
    d_{m}(z_2,...,z_d)\coloneqq\sum\limits_{\ell=1}^{\mu_{m}}    \frac{a_{\ell}^{(j_m)}}{(z_2 - b_{\ell 2}^{(j_m)}) ...(z_d - b_{\ell d}^{(j_m)})}.
\end{equation}
Applying a univariate rational reconstruction method from Section \ref{recuni} with the
 samples 
\begin{equation}\label{rat4}
r(k_1,\boldsymbol{0}_{d-1})=  \sum\limits_{m=1}^{M_1}  \frac{d_{m}(\boldsymbol{0}_{d-1})}{k_1 - b_{j_m 1}}, \quad k_1=-N,...,N,
\end{equation}
we reconstruct $M_1$ and the pairwise distinct poles $b_{j_m 1}$, $m=1,...,M_1$. Again, we have that $\mu_m \leq N$. Further, for a grid of points $(k_2,...,k_d) \in [-N,N]^{d-1} \cap \zz^{d-1}$, we compute the function values $d_{m}(k_2,...,k_d)$ as least squares solutions of systems of linear equations 
\begin{equation}\label{dl1}
r(k_1,k_2,...,k_d)=  \sum\limits_{m=1}^{M_1}  \frac{d_{m}(k_2,...,k_d)}{k_1 - b_{j_m 1}},  \quad k_1=-N,...,N.
\end{equation}
Thus, we obtain the function values $d_{m}(k_2,...,k_d)$ of     $(d-1)$-variate rational functions $d_{m}$, $m=1,...,M_1$,  on the grid $(k_2,...,k_d) \in [-N,N]^{d-1} \cap \zz^{d-1}$. That means that at the first step we reconstructed pairwise distinct poles $b_{j_m 1}$ of the 1st dimension and reduced the $d$-dimensional problem (\ref{rat}) to several (in total $M_1$ many) $(d-1)$-dimensional ones (\ref{dminus1r}), and using the values of the $d$-variate rational function (\ref{rat}) on the grid $[-N,N]^{d} \cap \zz^{d}$, we compute the values of the $(d-1)$-variate rational functions (\ref{dminus1r}) on the grid $[-N,N]^{d-1} \cap \zz^{d-1}$. 

Further, we proceed with recovering all poles with respect to 2nd
dimension, such that they create already correct pairs to the poles $b_{j_m 1}$. We apply similar ideas as above, but now instead of a function $r$ as in  (\ref{rat}) we work with $M_1$ functions $d_{m}$ as in (\ref{dminus1r}).
  Let us assume that for each $m=1,...,M_1$ we have $M_2^{(j_m)}$ pairwise distinct poles $b_{\ell 2}^{(j_m)}$. We denote them by $b_{\ell_{s_1}^{(1)} 2}^{(j_m)}$, $s_1=1,...,M_2^{(j_m)}$ and let  $\mu_{s_1}^{(j_m)}$  be their multiplicities. Then $\sum_{{s_1}=1}^{M_2^{(j_m)}} \mu_{s_1}^{(j_m)}= \mu_m$ and consequently $\sum_{m=1}^{M_1}  \sum_{{s_1}=1}^{M_2^{(j_m)}} \mu_{s_1}^{(j_m)} = \sum_{m=1}^{M_1}  \mu_m =M$ and we can rewrite rational functions $d_m$  in (\ref{dminus1r}) as
$$
    d_{m}(z_2,...,z_d) = \sum\limits_{s_1=1}^{M_{2}^{(j_m)}}     \frac{d_{s_1}^{(m)}(z_3,...,z_d)}{z_2-b_{\ell_{s_1}^{(1)} 2}^{(j_m)}},
    $$
where 
\begin{equation}\label{ds1}
d_{s_1}^{(m)}(z_3,...,z_d)\coloneqq \sum\limits_{s=1}^{\mu_{s_1}^{(j_m)}}  \frac{a_{s}^{(j_m,\ell_{s_1}^{(1)})}}{(z_3 - b_{s 3}^{(j_m,\ell_{s_1}^{(1)})}) ...(z_d - b_{s d}^{(j_m,\ell_{s_1}^{(1)})})}
\end{equation}
are $(d-2)$-variate rational functions and $a_{s}^{(j_m,\ell_{s_1}^{(1)})}$ are some coefficients. Applying again a univariate rational method  from  Section \ref{recuni} 
 $M_1\leq M$ times (for each $m=1,...,M_1$),  using the samples 
\begin{equation}\label{dm1}
    d_{m}(k_2,\boldsymbol{0}_{d-2})= \sum\limits_{s_1=1}^{M_{2}^{(j_m)}}     \frac{d_{s_1}^{(m)}(\boldsymbol{0}_{d-2})}{k_2-b_{\ell_{s_1}^{(1)} 2}^{(j_m)}}, \quad k_2=-N,...,N,
\end{equation}
we reconstruct  $M_{2}^{(j_m)}$ as well as the poles $b_{\ell_{s_1}^{(1)} 2}^{(j_m)}$, $s_1=1,...,M_2^{(j_m)}$, that are pairs to the poles $b_{j_m 1}$. Further,  we compute the function values $d_{s_1}^{(m)}(k_3,...,k_d)$ on the grid of points $(k_3,...,k_d) \in [-N,N]^{d-2} \cap \zz^{d-2}$ as least squares solutions of systems of linear equations 
\begin{equation}\label{dm1s}
 d_{m}(k_2, ...,k_d)=\sum\limits_{s_1=1}^{M_{2}^{(j_m)}}     \frac{d_{s_1}^{(m)}(k_3,...,k_d)}{k_2-b_{\ell_{s_1}^{(1)} 2}^{(j_m)}}, \quad k_2=-N,...,N,
\end{equation}
using the values $d_{m}(k_2, ...,k_d)$ on the grid $(k_2,...,k_d) \in [-N,N]^{d-1} \cap \zz^{d-1}$ computed at the previous step.

In the general case, at the step $p=3,...,d$, we want to reconstruct $M_p^{(j_m,\ell_{s_1}^{(1)},...,\ell_{s_{p-2}}^{(p-2)})}$ pairwise distinct poles $b_{\ell_{s_{p-1}}^{(p-1)},p}^{(j_m,\ell_{s_1}^{(1)},...,\ell_{s_{p-2}}^{(p-2)})}$, $s_{p-1}=1,...,M_p^{(j_m,\ell_{s_1}^{(1)},...,\ell_{s_{p-2}}^{(p-2)})}$ with respect to the dimension $p$ that are pairs to the poles $b_{\ell_{s_{p-2}}^{(p-2)},p-1}^{(j_m,\ell_{s_1}^{(1)},...,\ell_{s_{p-3}}^{(p-3)})}$. In order to do it, we consider the following representation of $(d-p+1)$-variate rational functions
$$
    d_{s_{p-2}}^{(m,s_1,...,s_{p-3})}(z_p,...,z_d)=\sum\limits_{s_{p-1}=1}^{M_p^{(j_m,\ell_{s_1}^{(1)},...,\ell_{s_{p-2}}^{(p-2)})}} \frac{ d_{s_{p-1}}^{(m,s_1,...,s_{p-2})}(z_{p+1},...,z_d)}{z_p-b_{\ell_{s_{p-1}}^{(p-1)},p}^{(j_m,\ell_{s_1}^{(1)},...,\ell_{s_{p-2}}^{(p-2)})} }, 
    $$
where $d_{s_{p-1}}^{(m,s_1,...,s_{p-2})}(z_{p+1},...,z_d)$ is a $(d-p)$-variate rational function of the form
\begin{equation}\label{dp1r}
d_{s_{p-1}}^{(m,s_1,...,s_{p-2})}(z_{p+1},...,z_d)\coloneqq \sum\limits_{s=1}^{\mu_{s_{p-1}}^{(j_m,\ell_{s_{1}}^{(1)},...,\ell_{s_{p-2}}^{(p-2)})}} \frac{a_{s}^{(j_m,\ell_{s_1}^{(1)},...,\ell_{s_{p-1}}^{(p-1)})}}{(z_{p+1} - b_{s (p+1)}^{(j_m,\ell_{s_1}^{(1)},...,\ell_{s_{p-1}}^{(p-1)})}) ...(z_d - b_{s d}^{(j_m,\ell_{s_1}^{(1)},...,\ell_{s_{p-1}}^{(p-1)})})},
\end{equation}
where $\mu_{s_{p-1}}^{(j_m,\ell_{s_{1}}^{(1)},...,\ell_{s_{p-2}}^{(p-2)})}$ describes the multiplicity of  $b_{\ell_{s_{p-1}}^{(p-1)},p}^{(j_m,\ell_{s_1}^{(1)},...,\ell_{s_{p-2}}^{(p-2)})}$, and $a_{s}^{(j_m,\ell_{s_1}^{(1)},...,\ell_{s_{p-1}}^{(p-1)})}$ are some coefficients (if $p=d$, $d_{s_{d-1}}^{(m,s_1,...,s_{d-2})}$ are just complex coefficients). By induction
\begin{align}
    \sum\limits_{m=1}^{M_1}  \sum\limits_{s_1=1}^{M_2^{(j_m)}} ... \sum\limits_{s_{p-2}=1}^{M_{p-1}^{(j_m,\ell_{s_1}^{(1)},...,\ell_{s_{p-3}}^{(p-3)})}} \left( \sum\limits_{s_{p-1}=1}^{M_p^{(j_m,\ell_{s_1}^{(1)},...,\ell_{s_{p-2}}^{(p-2)})}}  \mu_{s_{p-1}}^{(j_m,\ell_{s_{1}}^{(1)},...,\ell_{s_{p-2}}^{(p-2)})} \right)  =M. \label{mrec}
\end{align}
Applying one of the algorithms for rational approximation from Section \ref{recuni} $$\sum\limits_{s_{p-3}=1}^{M_{p-2}^{(j_m,\ell_{s_1}^{(1)},...,\ell_{s_{p-4}}^{(p-4)})}} M_{p-1}^{(j_m,\ell_{s_1}^{(1)},...,\ell_{s_{p-3}}^{(p-3)})}\leq M$$ (if $p=3$, $\sum_{m=1}^{M_1} M_2^{j_m}  \leq M$) times,  using the samples 
\begin{equation}\label{dsp2}
   d_{s_{p-2}}^{(m,s_1,...,s_{p-3})}(k_p,\boldsymbol{0}_{d-p})=\sum\limits_{s_{p-1}=1}^{M_p^{(j_m,\ell_{s_1}^{(1)},...,\ell_{s_{p-2}}^{(p-2)})}} \frac{ d_{s_{p-1}}^{(m,s_1,...,s_{p-2})}(\boldsymbol{0}_{d-p})}{k_p-b_{\ell_{s_{p-1}}^{(p-1)},p}^{(j_m,\ell_{s_1}^{(1)},...,\ell_{s_{p-2}}^{(p-2)})} }, \quad k_p=-N,...,N,
\end{equation}
we reconstruct  $M_p^{(j_m,\ell_{s_1}^{(1)},...,\ell_{s_{p-2}}^{(p-2)})}$ and the pairwise distinct poles $b_{\ell_{s_{p-1}}^{(p-1)},p}^{(j_m,\ell_{s_1}^{(1)},...,\ell_{s_{p-2}}^{(p-2)})} $, $s_{p-1}=1,...,M_p^{(j_m,\ell_{s_1}^{(1)},...,\ell_{s_{p-2}}^{(p-2)})}$. If $p\leq d-1$, we compute the values of rational functions (\ref{dp1r}) on the grid $(k_{p+1},...,k_d) \in [-N,N]^{d-p} \cap \zz^{d-p}$ as solutions to the following linear systems
\begin{equation}\label{dsp2s}
    d_{s_{p-2}}^{(m,s_1,...,s_{p-3})}(k_p,...,k_d)=\sum\limits_{s_{p-1}=1}^{M_p^{(j_m,\ell_{s_1}^{(1)},...,\ell_{s_{p-2}}^{(p-2)})}} \frac{ d_{s_{p-1}}^{(m,s_1,...,s_{p-2})}(k_{p+1},...,k_d)}{k_p-b_{\ell_{s_{p-1}}^{(p-1)},p}^{(j_m,\ell_{s_1}^{(1)},...,\ell_{s_{p-2}}^{(p-2)})} }, \quad k_p=-N,...,N,
\end{equation}
using the values $ d_{s_{p-2}}^{(m,s_1,...,s_{p-3})}(k_p,...,k_d)$ on the grid $(k_{p},...,k_d) \in [-N,N]^{d-p+1} \cap \zz^{d-p+1}$ computed at the previous step. Then the vectors of poles are reconstructed as
\begin{equation}\label{pol1}
\boldsymbol{b}_{m,s_1,...,s_{d-1}}=(  b_{j_m1},   b_{\ell_{s_1}^{(1)} 2}^{(j_m)},...,b_{\ell_{s_{d-1}}^{(d-1)}d}^{(j_m,\ell_{s_1}^{(1)},...,\ell_{s_{d-2}}^{(d-2)})} ), 
\end{equation}
and the frequency vectors  as $\boldsymbol{\lambda}_{m,s_1,...,s_{d-1}}=\frac{2\pi \mathrm{i}}{P} \boldsymbol{b}_{m,s_1,...,s_{d-1}}$
for  $m=1,...,M_1$, $s_1=1,...,M_2^{(j_m)}$, $s_{p-1}=1,...,M_p^{(j_m,\ell_{s_1}^{(1)},...,\ell_{s_{p-2}}^{(p-2)})}$, $p=3,...,d$. According to (\ref{mrec}) with $p=d$, the total number of the frequency vectors is exactly $M$. Note, that at the step $p=d$ we can finally reconstruct multiplicities $\mu_{s_{d-1}}^{(j_m,\ell_{s_{1}}^{(1)},...,\ell_{s_{d-2}}^{(d-2)})}$  of poles $b_{\ell_{s_{d-1}}^{(d-1)},d}^{(j_m,\ell_{s_1}^{(1)},...,\ell_{s_{d-2}}^{(d-2)})}$.

The coefficients $\gamma_j$, $j=1,...,M$ are reconstructed by (\ref{gamjd}), where $a_j$ are computed 
 as the solution to the following  linear system 
\begin{equation}\label{acoef}
    r(\boldsymbol{k})=\sum\limits_{j=1}^{M}      \frac{a_{j}}{(k_1 - b_{j 1}) ...(k_d - b_{j d})}, \quad \boldsymbol{k} \in [-N,N]^d.
\end{equation}


\begin{remark}\label{r51}
The concept of recovering poles and identifying correct pairs using the recursive dimension reduction approach can be represented as  \textit{"tree structure"}.  We have $M_1$ roots that correspond to the  $M_1$ pairwise distinct poles $b_{j_m 1}$, $m=1,...,M_1$ with respect to the first dimension.  Then each root $b_{j_m 1}$, $m=1,...,M_1$ has $M_2^{(j_m)}$ edges and children  nodes that correspond to the
$M_2^{(j_m)}$ pairwise distinct poles $b_{\ell_{s_1}^{(1)} 2}^{(j_m)}$, $s_1=1,...,M_2^{(j_m)}$ that are pairs to the $b_{j_m 1}$. That means that at the first level we have $\sum_{m=1}^{M_1} M_2^{(j_m)}$ poles with respect to the second dimension.  The index $^{(j_m)}$ in the notation $b_{\ell_{s_1}^{(1)} 2}^{(j_m)}$ shows to which root $b_{j_m 1}$ the poles $b_{\ell_{s_1}^{(1)} 2}^{(j_m)}$ belong. Then each node $b_{\ell_{s_1}^{(1)} 2}^{(j_m)}$, $s_1=1,...,M_2^{(j_m)}$, $m=1,...,M_1$  has $M_3^{(j_m,\ell_{s_1}^{(1)})}$ edges and children nodes $b_{\ell_{s_{2}}^{(2)},3}^{(j_m,\ell_{s_1}^{(1)})}$, $s_{2}=1,...,M_3^{(j_m,\ell_{s_1}^{(1)})}$. The notation $^{(j_m,\ell_{s_1}^{(1)})}$ means that a node $b_{\ell_{s_{2}}^{(2)},3}^{(j_m,\ell_{s_1}^{(1)})}$ is located at the 2nd level and its parent node and the root node are $b_{\ell_{s_1}^{(1)} 2}^{(j_m)}$ and   $b_{j_m 1}$, respectively.
In the general case, at the level $p-1$ for $p=4,...,d$ we have the nodes $b_{\ell_{s_{p-1}}^{(p-1)},p}^{(j_m,\ell_{s_1}^{(1)},...,\ell_{s_{p-2}}^{(p-2)})}$, $s_{p-1}=1,...,M_p^{(j_m,\ell_{s_1}^{(1)},...,\ell_{s_{p-2}}^{(p-2)})}$. The notation  $^{(j_m,\ell_{s_1}^{(1)},...,\ell_{s_{p-2}}^{(p-2)})}$ indicates all previous  parent nodes $b_{\ell_{s_{p-2}}^{(p-2)},p-1}^{(j_m,\ell_{s_1}^{(1)},...,\ell_{s_{p-3}}^{(p-3)})}$, ..., $b_{\ell_{s_1}^{(1)} 2}^{(j_m)}$ at levels $p-2,...,1$, respectively, and the root $b_{j_m 1}$. In Example \ref{ex62}, we visualize this tree structure.
\end{remark}



The main ideas of this method are presented in Algorithm \ref{alg4}. The algorithm  employs the full grid of the Fourier coefficients $ c_{\boldsymbol{k}}(f)$, $\boldsymbol{k} \in [-N,N]^d$ (i.e. $\mathcal{O}( N^d)$) and requires $\mathcal{O}(N M^2((d-1) M^2 + N^{d-1} ))$ flops. 
To reconstruct the poles $b_{j1}$, we apply the AAA method ones, which has the complexity $\mathcal{O}(N M^3)$. To determine the poles $b_{j p}$ for each $p=2,...,d$ we apply the AAA method at most $M$ times, which has the complexity $\mathcal{O}(N M^4)$. That means that we apply the AAA method in general at most $(d-1)M$ times, resulting in the computational cost of   $\mathcal{O}((d-1) N M^4)$. Finally, the complexity of the solutions of systems of linear equations (\ref{dl1}), (\ref{dm1s}) and (\ref{dsp2s}) is $\mathcal{O}(N^d M^2)$.

\section{Numerical experiments}
We apply now our algorithms for several examples. We have implemented all algorithms
with IEEE double precision arithmetic, using the free software Python. It will be convenient for us to say that the exponential sum (\ref{mul}) is given by the frequency matrix \(\bm \Lambda = (\lambda_{j\ell})_{j,l=1}^{M,d}\) and the coefficient vector \(\bm \gamma = (\gamma_j)_{j=1}^M\).  By  \(\tilde{\lambda}_{j\ell},\tilde \gamma\) and \(\tilde f\) we denote the reconstructed parameters \(\lambda_{j\ell},\gamma_j\) and the exponential sum \(f\), respectively. Then we define the relative reconstruction errors
\[ e(\bm \Lambda)\coloneqq  \max_{\ell=1,\dots,d} \frac{\displaystyle\max_{j=1,\dots, M}\abs{\lambda_{j\ell}-\tilde{\lambda}_{j\ell}}}{\displaystyle\max_{j=1,\dots,M}\abs\lambda_{j\ell}},\qquad e(\bm \gamma)\coloneqq \frac{\displaystyle\max_{j=1,\dots,M}\abs{\gamma_j-\tilde \gamma _j}}{\displaystyle \max_{j=1,\dots,M}\abs{\gamma_j}} \]
for  \(\bm \Lambda \) and  \(\bm \gamma \), respectively, and the relative  error for the exponential sums
\[e(f)\coloneqq \frac{\displaystyle \max_{\bm t}\abs{f(\bm t)-\tilde f(\bm t)}}{\displaystyle \max_{\bm t}\abs{f(\bm t)}}.\]
In order to estimate \(e(f)\) we choose \(51^d\) equidistant points \(\bm t\in[-10,10]^d\). In the following calculations we have determined all poles by solving the generalized eigenvalue problem \eqref{eig}. Although, using the Loewner pencil approach according to Section \ref{mpa}, we have observed equally good results.

\begin{figure}[h!]
  \centering
    \includegraphics[width=0.35\linewidth]{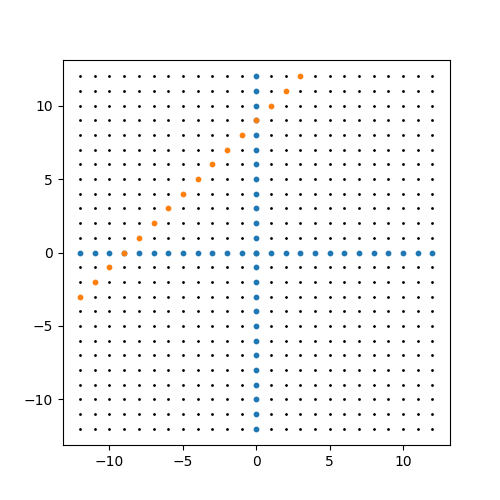}
  \caption{ The sparse grid of indices of the Fourier coefficients used by Algorithm \ref{alg1} for the reconstruction of $f_1$ from  Example \ref{ex61}. The blue dots visualize the coefficients, which are used for separate computation of the frequencies, while the orange dots indicate the coefficients used for determining the correct pairing. The grey dots represent all coefficients from the 2D cube $[-N,N]^2$.  }
  \label{2dex}
\end{figure}

\begin{figure}[h!]
  \centering
    \includegraphics[width=0.5\linewidth]{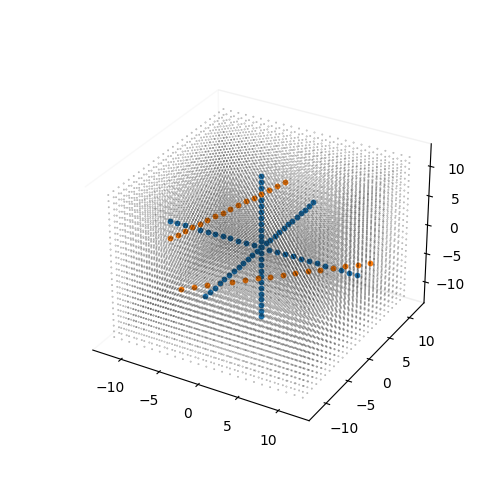}
  \caption{The sparse grid of indices of the Fourier coefficients used by Algorithm \ref{alg3} for the reconstruction of  $f_2$ from  Example \ref{ex61}. The blue dots visualize the coefficients, which are used for separate computation of the frequencies, while the orange dots indicate the coefficients used for determining the correct pairing. The grey dots represent all coefficients from the 3D cube $[-N,N]^3$.}
  \label{2dex1}
\end{figure}

\begin{example}
\label{ex61}
Let us for now confine to the setting \(\lambda_{j\ell}\neq \lambda_{i\ell}\) for \(j\neq i\)  and \(\ell=1,\dots,d\). We consider a bivariate exponential sum \(f_1\) of order 5 defined by $\bm{\gamma}_1$ and $\Lambda_1$ and a 3-variate exponential sum  \(f_2\) of order 5  given through $\bm{\gamma}_2$ and $\Lambda_2$:
\[\resizebox{\textwidth}{!}{%
$\bm{\gamma}_1 = \begin{pmatrix}
3\\2\\1\\2\\1
\end{pmatrix},\quad \bm \Lambda_1 =\begin{pmatrix}
\sqrt{2.21}i&3.33i\\
-5.63i&-\sqrt 5i\\
-3.47i&\sqrt{6}i\\
-\sqrt{7.1}i&-4.5i\\
0.46i&-9.44i	
\end{pmatrix},\quad \bm \gamma_2=\begin{pmatrix}
-1\\-2\\-3\\1\\2\\3
\end{pmatrix},\quad \bm \Lambda_2 = \begin{pmatrix}
-2-3i&\sqrt{\pi}(-1+i)&0.5i\\
-1+\sqrt{20}i&-3+i&-1+i\\
3i&-4+0.5i&1.22i\\
-2+3i&\sqrt{\pi}(-1-i)&-0.5i\\
-1-\sqrt{20}i&-3-i&-1-i\\
-3i&-4-0.5i&-1.22i\\
\end{pmatrix}.
$}\]
We present the reconstruction results by Algorithms \ref{alg3} (for 2D case we apply Algorithm \ref{alg1}) and \ref{alg4} as well as the used parameters $N$, $P$ and $\tau$ in   Table  \ref{tab1}.  The corresponding sparse grids used by Algorithms \ref{alg1} and \ref{alg3} are presented in Figures \ref{2dex} and \ref{2dex1}, respectively.
\begin{table}[h!]
\begin{center}
\caption{\small Reconstruction errors for Example \ref{ex61}}
 \label{tab1}
\begin{tabular}{|l|c|c|c||c|c|c||c|c|c|}\cline{5-10}
\multicolumn{4}{l}{ }&\multicolumn{3}{|c||}{Algorithm \ref{alg3}}&\multicolumn{3}{|c|}{Algorithm \ref{alg4}}\\\hline
&P&N&\(\tau\)&\(e(\bm \Lambda)\)&\(e(\bm\gamma)\)&\(e(f)\)&\(e(\bm \Lambda)\)&\(e(\bm\gamma)\)&\(e(f)\)\\\hline
$f_1$ &4&15&7&8.182e-14&3.212e-13&8349e-13&8.182e-14&3.015e-13&7.464e-13\\
$f_2$ &5&15&4&5.202e-10&5.344e-10&2.116e-09&
1.401e-13&8.134e-14&1.670e-12\\\hline
\end{tabular}
\end{center} 
  \end{table}
\end{example}

\begin{example}
\label{ex62}
We turn now to the general case, allowing repetitive frequencies. In this case we apply only Algorithm \ref{alg4}. Define a 4-variate exponential sum $f_3$ of order 9 by
$$
\bm{\Lambda}_3=\begin{pmatrix}
2+2i&0.2i&i&1\\
2+2i&0.2i&i&-1\\
2+2i&-2&1+i&i\\
2+2i&-2&1+i&-2i\\
2+2i&-2&1+i&3i\\
3+i&-\pi&-3&-\sqrt\pi i\\
3+i&-\pi&1&2i\\
3+i&-\pi&1&-4\\
3+i&0.2i&1+i&\sqrt{20}i
\end{pmatrix}.
$$
 and \(\gamma_3=(1,...,1)^T \in \cc^9 \).  Further, let $f_4$ be a 3-variate exponential sum of order 4 given by
\[ \gamma_4 =\begin{pmatrix} 1\\-4\\-2\\2
\end{pmatrix},\quad \bm\Lambda_4=\begin{pmatrix}
-1.47-0.27i&-1.87-0.57i&-1.35+4.61i\\
-1.47-0.27i&-1.87-0.57i&-1.26-2.58i\\
-1.47-0.27i&-0.84+7.53i&-1.75-1.33i\\
-0.60+4.86i&-0.13+5.05i&-0.12+8.34i\\
\end{pmatrix}.\]
Then the used parameters as well as the relative reconstruction errors by Algorithm \ref{alg4} can be found in Table  \ref{tab2}.
\begin{table}[h!]
\begin{center}
 \caption{\small Reconstruction errors for Example \ref{ex62}}
 \label{tab2}
\begin{tabular}{|c|c|c||c|c|c|}\hline
&\(P\)&\(N\)&\(e(\bm \Lambda)\)&\(e(\gamma)\)&\(e(f)\)\\\hline
$f_3$&2.4&10& 9.0480e-13&1.0865e-12&4.0311e-11\\
$f_4$&1&10&1.6710e-15&1.0215e-15&4.9217e-14\\\hline
\end{tabular}
\end{center}
  \end{table}
The recursive routine for recovering the frequencies of the exponential sum $f_3$
using Algorithm \ref{alg4}
follows the tree pattern
shown below:

\begin{center}
    \begin{tikzpicture}
		\node () at (-1.5,-1){roots};
		\node () at (-1.5,-2){level 1};
		\node () at (-1.5,-3){level 2};
		\node () at (-1.5,-4){level 3};        
        
        \node (11) at (5.25/2,-1) {\(2+2i\)};
        \node (12) at (10.3125,-1) {\(3+i\)};

        \node (21) at (0.75,-2) {\(0.2i\)};
        \node (22) at (4.5,-2) {\(-2\)};
        \node (23) at (17.25/2,-2) {\(-\pi\)};
        \node (24) at (12,-2) {\(0.2i\)};
        
        \node (31) at (0.75,-3) {\(i\)};
        \node (32) at (4.5,-3) {\(1+i\)};
        \node (33) at (7.5,-3) {\(-3\)};
        \node (34) at (9.75,-3) {\(1\)};
        \node (35) at (12,-3) {\(1+i\)};

        \node (41) at (0,-4) {\(1\)};
        \node (42) at (1.5,-4) {\(-1\)};
        \node (43) at (3,-4) {\(i\)};
        \node (44) at (4.5,-4) {\(-2i\)};
        \node (45) at (6,-4) {\(3i\)};
		\node (46) at (7.5,-4) {\(-\sqrt \pi i\)};
		\node (47) at (9,-4) {\(2i\)};
		\node (48) at (10.5,-4) {\(-4\)};				
		\node (49) at (12,-4) {\(\sqrt{20}i\)};
        \draw (31) -- (41);
        \draw (31) -- (42);
        \draw (32) -- (43);
        \draw (32) -- (44);
        \draw (32) -- (45);
        \draw (33) -- (46);
        \draw (34) -- (47);
        \draw (34) -- (48);
        \draw (35) -- (49);

        \draw (21) -- (31);
        \draw (22) -- (32);
        \draw (23) -- (33);
        \draw (23) -- (34);
        \draw (24) -- (35);
        
        \draw (11) -- (21);
        \draw (11) -- (22);
        \draw (12) -- (23);
        \draw (12) -- (24);
        
    \end{tikzpicture}
    
\end{center}
\end{example}
\begin{example} 
\label{ex63}
We choose the same examples as in \cite{PT13mult}, i.e. we define a bivariate exponential sum $f_5$, a 3-variate exponential sum  $f_6$, and 4-variate exponential sum $f_7$ of orders 8 by
\[
 \bm \gamma_5 =\bm \gamma_6 =\bm \gamma_7 =\begin{pmatrix} 1+i\\2+3i\\5-6i\\0.2-i\\1+i\\2+3i\\5-6i\\0.2-i\end{pmatrix},
\quad
\bm\Lambda_5 = \begin{pmatrix}
0.1i&1.2i\\
0.19i&1.3i\\
0.3i&1.5i\\
0.35i&0.3i\\
-0.1i&1.2i\\
-0.19i&0.35i\\
-0.3i&-1.5i\\
-0.3i&0.3i
\end{pmatrix},
\]
\[\bm \Lambda_6=
\begin{pmatrix}
0.1i&1.2i&0.1i\\
0.19i&1.3i&0.2i\\
0.4i&1.5i&1.5i\\
0.45i&0.3i&-0.3i\\
-0.1i&1.2i&0.1i\\
-0.19i&0.35i&-0.5i\\
-0.4i&-1.5i&0.25i\\
-0.4i&0.3i&-0.3i
\end{pmatrix},\quad \bm \Lambda_7 = \begin{pmatrix}
0.1i&1.2i&0.1i&0.45i\\
0.19i&1.3i&0.2i&1.5i\\
0.3i&1.5i&1.5i&-1.3i\\
0.45i&0.3i&-0.3i&0.4i\\
-0.1i&1.2i&0.1i&-1.5i\\
-0.19i&0.35i&-0.5i&-0.45i\\
-0.4i&-1.5i&0.25i&1.3i\\
-0.4i&0.3i&-0.3i&0.4i
\end{pmatrix}. \]
The corresponding recovery results by Algorithm \ref{alg4} are presented in Table \ref{tab3}.
\begin{table}[h!]
\begin{center}
 \caption{\small Reconstruction errors for Example \ref{ex63}}
\begin{tabular}{|c|c|c||c|c|c|}\hline
 \label{tab3}
&P&N&\(e(\bm\Lambda)\)&\(e(\bm\gamma)\)&\(e(f)\)\\\hline
$f_5$&60&15&1.2881e-14&4.1350e-14&3.5417e-14\\
$f_6$&60&15&2.5387e-15&3.1120e-14&1.8197e-14\\
$f_7$&60&15&1.5535e-14&6.6704e-14&8.3695e-14\\\hline
\end{tabular}
\end{center}
  \end{table}

\end{example}  

\section{Conclusions}
We develop two new methods for the recovery of multivariate exponential sums,  using a finite set of their Fourier coefficients as input data. Exploiting the rational structure of the Fourier coefficients, we replace the multivariate exponential recovery problem by the multivariate rational interpolation problem. Further, this special multivariate rational interpolation problem can be solved applying the sparse grid approach or the recursive dimension reduction method.

As input information, instead of the Fourier coefficients $ c_{\boldsymbol{k}}(f)$ of $f$ in (\ref{mul}), we can also use its function values $f(\boldsymbol{\ell})$, $\boldsymbol{\ell} \in [0,2N-1]^d \cap \mathbb{Z}^d$. 
 In  \cite{ DPR23,  DPP21},  we  developed our univariate  method \textit{ESPIRA}  (Estimation of Signal Parameters via Iterative Rational Approximation) by replacing Fourier coefficients of $f$ in (\ref{unexp}) by the corresponding DFT (discrete Fourier Transform) and showing that modified DFT of univariate exponential sums also have a rational structure. 
Note that a similar idea can be also applied in the multivariate setting. Indeed, having the function values $f(\boldsymbol{\ell})$, $\boldsymbol{\ell} \in [0,2N-1]^d \cap \mathbb{Z}^d$,  we can compute their $d$-variate DFT  $\hat{f}(\boldsymbol{k})$, $\boldsymbol{k} \in [0,2N-1]^d \cap \mathbb{Z}^d$ by
\begin{align*}
    \hat{f}(\boldsymbol{k}) &= \sum \limits_{\boldsymbol{\ell} \in [0,2N-1]^d \cap \mathbb{Z}^d} f(\boldsymbol{\ell}) \mathrm{e}^{- \frac{2 \pi \mathrm{i}}{2N} \langle \boldsymbol{\ell}, \boldsymbol{k} \rangle}  = \sum\limits_{j=1}^{M} \gamma_j  \sum \limits_{\boldsymbol{\ell} \in [0,2N-1]^d \cap \mathbb{Z}^d} \mathrm{e}^{\langle \boldsymbol{\lambda}_j, \boldsymbol{\ell} \rangle} \mathrm{e}^{- \frac{2 \pi \mathrm{i}}{2N} \langle \boldsymbol{\ell}, \boldsymbol{k} \rangle}  \\
    & = \sum\limits_{j=1}^{M} \gamma_j \prod\limits_{m=1}^d   \sum\limits_{\ell_m=0}^{2N-1} \mathrm{e}^{(\lambda_{j m}-\frac{2 \pi \mathrm{i}}{2N}  k_1)\ell_m }    = \mathrm{e}^{\frac{2 \pi \mathrm{i}}{2N} \langle \boldsymbol{1}, \boldsymbol{k} \rangle}  \sum\limits_{j=1}^{M} \gamma_j \prod \limits_{m=1}^d   \frac{ 1-\mathrm{e}^{\lambda_{j m} 2N} }{  \mathrm{e}^{\frac{2 \pi \mathrm{i}}{2N} k_m}-\mathrm{e}^{\lambda_{j m} }},
\end{align*}
where $\boldsymbol{1}=(1,...,1) \in  \cc^d$ and $\lambda_{j m} \neq \frac{2 \pi \mathrm{i}}{2N} k_m$.
From the last formula we get
$$
    \mathrm{e}^{-\frac{2 \pi \mathrm{i}}{2N} \langle \boldsymbol{1}, \boldsymbol{k} \rangle}  \hat{f}(\boldsymbol{k})  =  \sum\limits_{j=1}^{M} \gamma_j   \prod \limits_{m=1}^d   \frac{ 1-\mathrm{e}^{\lambda_{j m} 2N} }{  \mathrm{e}^{\frac{2 \pi \mathrm{i}}{2N} k_m}-\mathrm{e}^{\lambda_{j m} }}.
$$
Defining a $d$-variate rational function $r(\boldsymbol{z})$ as in (\ref{rat}) with $b_{j m}= \mathrm{e}^{\lambda_{j m} }$ and $a_j = \gamma_j  \prod_{m=1}^d (1-\mathrm{e}^{\lambda_{j m} 2N})$, 
we get that $r$ satisfies the following interpolation condition 
\begin{equation}\label{dftint}
r(\mathrm{e}^{\frac{2 \pi \mathrm{i}}{2N} k_1}, ..., \mathrm{e}^{\frac{2 \pi \mathrm{i}}{2N} k_d}) = \mathrm{e}^{-\frac{2 \pi \mathrm{i}}{2N} \langle \boldsymbol{1}, \boldsymbol{k} \rangle}  \hat{f}(\boldsymbol{k}), \quad \boldsymbol{k} \in [0,2N-1]^d \cap \mathbb{Z}^d.
\end{equation}
Thus, we can reformulate the multivariate exponential recovery problem as  multivariate rational interpolation problem (\ref{dftint}) and apply techniques developed in Section \ref{spd}. More detailed investigation of this question as well as the recovery from given noisy data (Fourier coefficients $c_{\boldsymbol{k}}(f)(1 + \varepsilon_{\boldsymbol{k}})$
or function values $f({\boldsymbol{k}})(1 + \varepsilon_{\boldsymbol{k}})$, where $\varepsilon_{\boldsymbol{k}}$ is a vector which components are either i.i.d. random variables drawn from a standard normal
distribution with mean value 0 or from a uniform distribution with mean value 0)
 we leave for further research.


\section*{Statements and Declarations}
\textbf{Competing Interests:} The authors declare that they have no competing interests.

\small

\begin{thebibliography}{99999}

\bibitem{sturm2023}
R. Ait El Manssour,  M. Härkönen, B. Sturmfels:  Linear PDE with constant coefficients,  Glasgow Mathematical Journal,   \textbf{65}(S1),  2023, S2–S27.  doi:10.1017/S0017089521000355


\bibitem{ALI17}
A. C. Antoulas, S. Lefteriu,  A. C. Ionita, A tutorial introduction to the Loewner framework for model reduction,  in Model Reduction and Approximation,  SIAM,  Philadelphia,  2017,   335--376.  doi:10.1137/1.9781611974829.ch8




\bibitem{CL2018}
A.  Cuyt,  W.-s.  Lee,   Multivariate exponential analysis from the minimal
number of samples,   Adv.  Comput.  Math.,  \textbf{44},  2018, 987–1002.
https://doi.org/10.1007/s10444-017-9570-8


\bibitem{DB13}
A.  Damle,  G.  Beylkin, T.  Haut,  L.  Monzón,
Near optimal rational approximations of large data sets,
Appl.  Comput. Harmon. Anal.,
\textbf{ 35}(2), 
2013,
251-263.
https://doi.org/10.1016/j.acha.2012.08.011



\bibitem{DPP21}
N.~Derevianko,  G.~Plonka, and M.~Petz,  
\newblock From ESPRIT to ESPIRA: Estimation of signal parameters by iterative rational approximation, {\em IMA J. Numer. Anal.} \textbf{43}(2), 2023, 789--827. https://doi.org/10.1093/imanum/drab108

\bibitem{DPR23}
N.~Derevianko,  G.~Plonka,  and R. Razavi,
\newblock ESPRIT versus ESPIRA for reconstruction of short cosine sums and its application, {\em Numer. Algor.},  \textbf{92}, 2023,  437--470.


\bibitem{DP21}
N.~Derevianko,  G.~Plonka,  
Exact reconstruction of extended exponential sums using rational approximation of their Fourier coefficients,  Anal.  Appl.,  \textbf{20} (3),  2021,  543--577.   doi: 10.1142/S0219530521500196


\bibitem{D24}
N.~Derevianko,  
Recovery of rational functions via Hankel pencil method and sensitivities of the poles, Arxiv preprint: https://arxiv.org/abs/2406.13192


\bibitem{DI17}
B. Diederichs and A. Iske, Projection-based multivariate frequency estimation, 2017 International Conference on Sampling Theory and Applications (SampTA), Tallinn, Estonia, 2017, pp. 360-363, doi: 10.1109/SAMPTA.2017.8024446.









\bibitem{G2022}
S.  Jiang,   L.  Greengard,  Approximating the Gaussian as a sum of exponentials and its spplications to the fast Gauss transform, 
Commun.  Comput.  Phys., \textbf{31},   2022,  1-26.
doi:10.4208/cicp.OA-2021-0031



\bibitem{B23}
R. Katz, N. Diab, D. Batenkov, Decimated Prony's Method for Stable Super-Resolution, IEEE Signal Processing Letters, \textbf{30}, 2023, 1467-1471. 10.1109/LSP.2023.3324553

\bibitem{KGA}
D. S. Karachalios, I. V. Gosea, A. C. Antoulas, The Loewner framework for system identification and reduction. In P. Benner, S. Grivet-Talocia, A. Quarteroni, G. Rozza, W. H. A. Schilders, and L. M. Silveira, editors, Handbook on Model Reduction, volume I of Methods and Algorithms in press (2020).



\bibitem{KP2016}
S.  Kunis, T.  Peter, T.  Römer,  U.  von der Ohe,
A multivariate generalization of Prony's method,
Linear Algebra Appl.,
\textbf{ 490},
2016,  31-47.
https://doi.org/10.1016/j.laa.2015.10.023

\bibitem{L20}
W. Li, W. Liao, and A. Fannjiang, Super-Resolution Limit of the ESPRIT Algorithm, IEEE Transactions on Information Theory, \textbf{ 66}(7), 2020,  4593 - 4608. 10.1109/TIT.2020.2974174


\bibitem{AAA}
Y. Nakatsukasa, O. Sete, and L.N. Trefethen,  The AAA algorithm for rational approximation.
SIAM J. Sci. Comput., \textbf{40}(3), A1494–A1522, 2018.

\bibitem{M25}
H. N. Mhaskar,, S. Kitimoon† and Raghu G. Raj, Robust and Tractable Multidimensional Exponential Analysis, 2025. https://arxiv.org/pdf/2404.11004

\bibitem{PLH04}
J.-M.  Papy,  L.  De Lathauwer,  S.  Van Huffel, Exponential data fitting using multilinear algebra: the decimative case,  J. Chemometrics, \textbf{23}, 2009,  341-351.  https://doi.org/10.1002/cem.1212





\bibitem{PT13mult}
D.  Potts,  M.  Tasche,  Parameter estimation for multivariate exponential sums,  Electron. Trans. Numer.
Anal.,  \textbf{40},  2013,  204–224.  https://etna.math.kent.edu/volumes/2011-2020/vol40/abstract.php?vol=40\&pages=204-224



















































\bibitem{PT13}
 D. Potts, M. Tasche, Parameter estimation for nonincreasing exponential sums by Prony-like methods, Linear Algebra Appl. \textbf{439}, 2013, 1024--1039.
 

\bibitem{PV2020}
J.  Prestin,  H.  Veselovska,  Prony-type polynomials and their common zeros,  Front.  Appl.  Math.  Stat.,   \textbf{6}:16,  2020,   21 pp.   doi: 10.3389/fams.2020.00016


\bibitem{R22}
A. C. Rodriguez, L. Balicki, and S. Gugercin, Serkan, The p-AAA Algorithm for Data-Driven Modeling of Parametric Dynamical Systems, SIAM Journal on Scientific Computing,  \textbf{45}, 2023, A1332-A1358. https://doi.org/10.1137/20M1322698

\bibitem{RT89}
R. Roy, T. Kailath,  ESPRIT – estimation of signal parameters via rotational invariance techniques, IEEE Trans.  Acoust.  Speech
Signal Process.,   \textbf{37}(7),  1989, 984–995. 10.1109/29.32276












\bibitem{S2018}
T.  Sauer,  Prony’s method in several variables: symbolic solutions by universal interpolation,  J. Symb.
Comput.,  \textbf{84,}  2018,  95–112.  https://doi.org/10.1016/j.jsc.2017.03.006

\bibitem{SD2007}
P. Shukla,  P. L. Dragotti,  Sampling schemes for multidimensional signals with finite rate of
innovation,  IEEE Trans.  Signal Process.,   \textbf{55}(7),  2007,  3670–3686.  doi: 10.1109/TSP.2007.894259

 \bibitem{WDT21}
H.~Wilber, A.~Damle, and A.~Townsend,  Data-driven algorithms for signal processing with trigonometric rational functions, 
\newblock {\em SIAM J. Sci. Comput.}, \textbf{44}(3), 2022, C185--C209.

\bibitem{YFS2006}
N.  Yilmazer,  R.  Fernandez-Recio,  T.K.  Sarkar,  Matrix pencil method for simultaneously estimating
azimuth and elevation angles of arrival along with the frequency of the incoming signals.  Digital
Signal Process.,   \textbf{16}(6),  2006,  796–816.  https://doi.org/10.1016/j.dsp.2006.05.009.

\bibitem{Y22}  L.  Ying,  Pole recovery from noisy data on imaginary axis, J. Sci. Comput.,  \textbf{92}(107),  2022.  https://doi.org/10.1007/s10915-022-01963-z










 
 
 
 
\end{thebibliography}

\end{document}